\documentclass[journal,onecolumn]{IEEEtran}
\makeatletter
\def\ps@headings{%
\def\@oddhead{\mbox{}\scriptsize\rightmark \hfil \thepage}%
\def\@evenhead{\scriptsize\thepage \hfil\leftmark\mbox{}}%
\def\@oddfoot{}%
\def\@evenfoot{}}
\makeatother
\usepackage{latexsym}
\usepackage{amssymb}
\usepackage{epsfig}
\usepackage{amsthm} 
\usepackage{subfigure}
\usepackage{graphics,color}
\usepackage{graphicx}
\usepackage{amsmath}
\usepackage{algorithm}
\usepackage{algorithmic}
\usepackage{cite}
\usepackage{comment}
\usepackage{url}
\usepackage{xspace}

\newtheorem{theorem}{Theorem}
\newtheorem{lemma}{Lemma}

\newtheorem{corollary}{Corollary}

\theoremstyle{definition}
\newtheorem{definition} {Definition}
\newtheorem{remarks} {Remark}
\newtheorem{example}{Example}
\newtheorem{claim}{Claim}
\newtheorem{fact}{Fact}

\renewcommand{\qed}{\nobreak \ifvmode \relax \else
      \ifdim\lastskip<1.5em \hskip-\lastskip
      \hskip1.5em plus0em minus0.5em \fi \nobreak
      \vrule height0.2em width0.5em depth0.4em\fi}
\IEEEoverridecommandlockouts
\begin{document}

\title{A General Class of Throughput Optimal Routing Policies in Multi-hop Wireless Networks}
\author{M. Naghshvar, H. Zhuang, and T. Javidi \\
Department of Electrical and Computer Engineering, \\
University of California San Diego, \\
La Jolla, CA 92093 USA, \\
\{naghshvar, hzhuang, tjavidi\}@ucsd.edu
}

\maketitle

\begin{abstract}
This paper considers the problem of throughput optimal routing/scheduling in a multi-hop constrained queueing network with random connectivity whose special case includes opportunistic multi-hop wireless networks and input-queued switch fabrics. 
The main challenge in the design of throughput optimal routing policies is closely related to identifying appropriate and universal Lyapunov functions with negative expected drift. The few well-known throughput optimal policies in the literature are constructed using simple quadratic or exponential Lyapunov functions of the queue backlogs and as such they seek to balance the queue backlogs across network independent of the topology.

By considering a class of continuous, differentiable, and piece-wise quadratic Lyapunov functions, this paper provides a large class of throughput optimal routing policies.
The proposed class of Lyapunov functions allow for the routing policy to control the traffic along short paths for a large portion of state-space while ensuring a negative expected drift. This structure enables the design of a large class of routing policies.
In particular, and in addition to recovering the throughput optimality of the well known backpressure routing policy, an opportunistic routing policy with congestion diversity is proved to be throughput optimal.

\end{abstract}



\section{Introduction}

\thispagestyle{empty}

This paper considers the problem of throughput optimal routing/scheduling in a general constrained queueing network with random connectivity
whose special case includes opportunistic routing in multi-hop wireless network and input-queued switch scheduling.
While it is often possible to intuitively design and propose various routing/scheduling policies, providing theoretical guarantees for the corresponding controlled Markov chains is far from straight forwards with the exception of the throughput optimality of backpressure routing \cite{TassEph92} and maximum weight scheduling \cite{McKeown99}.
These guarantees are obtained using Foster-Lyapunov Theorem which ensures the stability of a controlled Markov chain if a Lyapunov function with negative expected drift is shown to exist.
More specifically, the throughput optimal backpressure-based policies \cite{TassEph92,Neely09,Xi06,sarkar08,Ying08} 
as well as maximum weight schedules \cite{Tassiulas93,McKeown99} 
are reverse-engineered to be the very rule under which the known quadratic Lyapunov function is ensured a negative expected drift.

\begin{figure}[htp]
\centering
\includegraphics[width=0.6\textwidth]{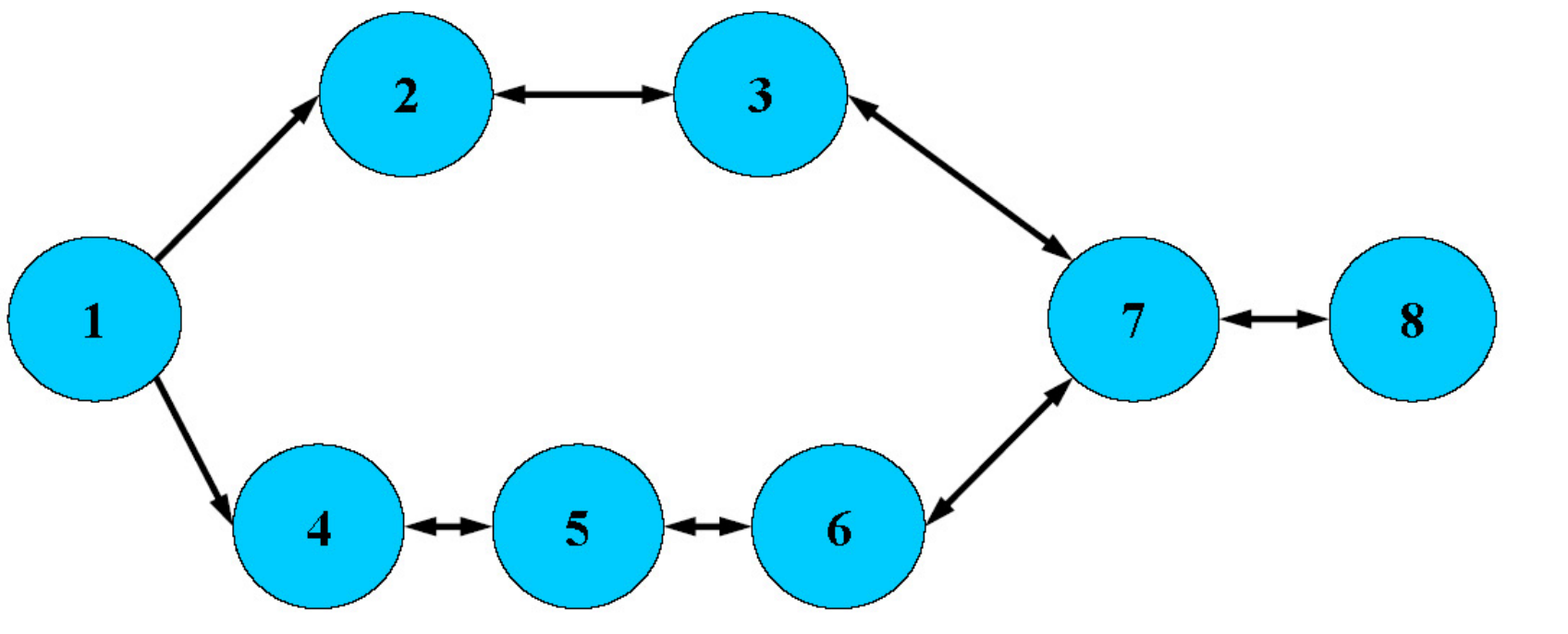}    
\caption{A network of eight nodes. Packets are to be routed from node~1 to node~8.}
\label{fig:NetIntro}
\end{figure}

While reverse engineering routing/scheduling in this function has the advantage of obtaining theoretical guarantees, it may result in schemes with undesirable structure. 
In particular, under the strict Schur-convexity of quadratic Lyapunov function \cite{TassEph92,McKeown99} (as well as the exponential Lyapunov functions \cite{Hassibi08})
with respect to the (weighted) backlog vector,
the negativity of the expected drift is only achieved when nodes with large queues are prioritized in favor of those with small number of buffered packets (e.g. a node with small backlog must refrain from routing packets to a neighbor with large backlog).
This very need to ensure a negative drift of the Lyapunov function (equivalently to balance the queues in a network), goes against the intuition behind many promising routing/scheduling schemes.
For instance, consider the wired network in Fig.~\ref{fig:NetIntro} where packets are to be routed from node~1 to node~8. 
It is intuitively desirable for the routing decisions in this network to be 
such that the bottle-neck link (7,8) is maximally utilized.
Indeed, in Subsection~\ref{backORCD} we discuss an opportunistic routing policy (ORCD) which attempts to achieve this goal.
However, these very intuitive properties cause a positive expected drift in the quadratic Lyapunov function in an infinite number of states.
This means that theoretical guarantee for this algorithm requires a significantly different approach (non-Schur-convex Lyapunov function).

In this paper, we provide a large class of throughput optimal policies by considering a class of piece-wise quadratic Lyapunov functions. 
The proposed class of Lyapunov functions are constructed by grouping the queues based on their relative size and the network topology and as such are not strictly Schur-convex.
This allows for the Lyapunov function to have an expected negative drift even when packets are routed from a node with small backlog to one with large backlog so long as the queues are grouped together. 
We will see that the proposed class of Lyapunov functions establish 
the throughput optimality of a large class of routing/scheduling policies by indirectly incorporating the critical information about topology.  
In particular, we specialize our result to recover the throughput optimality of two known routing 
policies, backpressure (already known to be throughput optimal) and ORCD (discussed above and whose throughput optimality 
only was conjectured in \cite{parul07}). 

Before we close, we note that using the methodology in this paper, it is always possible to find uncountably 
many throughput optimal routing/scheduling policies among which most will 
suffer from an unreasonable complexity and overhead. In light of this observation, 
we believe that (even though beyond the scope of this paper) 
the contribution and utility of the proposed Lyapunov construction is of two folds: 1) The proposed class of Lyapunov functions can be used systematically to establish the throughput optimality of many 
intuitive solutions which do not locally balance the backlogs in the network. Furthermore, 2) the nature of the constructed Lyapunov function concretely establishes the intuition that queue stability in a network is \emph{only} affected by the control applied at the boundaries of the state space where inevitable idling is likely to occur (these are the states in which one or more of the queues are near empty while others have extremely large backlogs). The consequence of the former is more flexibility when proving the throughput optimality of various routing/scheduling solutions such as the multi-hop schemes that would favor short paths, while the latter establishes a qualitative characterization of throughput optimality.  


The remainder of this paper is organized as follows.
In Section~\ref{ProblemFormulation}, we formulate the problem for the general case of routing/scheduling in constrained queueing network with multiple destinations.
For ease of exposition, the results are first presented in Section~\ref{SingleDst} for the multi-hop routing in a single destination network 
with orthogonal channels.
The extensions of the results to the general constrained queueing is provided in Section~\ref{multidst} where we also show that
scheduling for input-queued switches is a special case of our framework.
Section~\ref{Design} discusses the structure of our proposed piece-wise quadratic Lyapunov function and provides a (alternative) proof of throughput optimality for some of the existing routing/scheduling policies.  
Finally, we conclude the paper and discuss future work in Section~\ref{Discussion}.

We close this section with a note on the notations used.
Let $[ x ] ^{+} = \max \{x,0 \}$. The indicator function $\mathbf{1}_{\{ X \}}$ takes the value $1$ whenever event $X$ occurs, and $0$ otherwise. For any set $S$, $\left| S \right|$ denotes the cardinality of $S$, while for any vector $\boldsymbol{v}$, $\left\| \boldsymbol{v} \right\|$ denotes the euclidean norm of $\boldsymbol{v}$. For matrices 
$A$ and $B$, let $\langle A, B\rangle =\sum_{i,j} A_{ij}B_{ij}$ denote the inner product. 
For any set $S$, $int(S)$ is the set of all interior points of $S$.
When dealing with a sequence of sets $C_1,C_2,\dots$, we define $C^i=\cup_{j=1}^{i} C_j$.
Lastly, we use bold letters to discriminate vectors from scalar quantities as well as their components.


\section{Problem Formulation}
\label{ProblemFormulation}


We consider a time slotted system with slots indexed by $t \in \{ 0,1,2,\ldots \}$ where slot $t$ refers to the time interval $[t,t+1)$. There are $N$ nodes in the network labeled by $\Omega:=\{ 1,2,\ldots,N \}$. 
We denote the set of all destinations by $\mathcal{D}$, $\mathcal{D} \subseteq \Omega$.

Let random variable $A_i^d(t)$ represent the amount of data (in units of packets) that exogenously arrives to node $i$ and destined for node $d$, $d \in \mathcal{D}$ during time slot $t$. Arrivals are assumed to be i.i.d. over time and bounded by a constant $A_{\max}$. 
All packets destined for node $d$, $d \in \mathcal{D}$ are referred to as \emph{commodity $d$ packets}.
Let $\lambda_i^d = \mathbb{E} [A_i^d(t)]$ denote the exogenous arrival rate of commodity $d$ packets to node $i$. 
We define $\boldsymbol{\lambda}=[\lambda_i^d]_{d \in \mathcal{D}, i \in \Omega}$ to be the arrival rate matrix (of size $|\mathcal{D}| \times N$).  
We assume that each node maintains a separate buffer (with infinite queuing space) for each destination
in which packets that arrive exogenously at that node as well as packets routed to that node from other nodes in the network are queued. 
Without loss of generality, we assume that after a packet is successfully received at its destination, the packet would be ejected from the network. 
Let $Q_i^d(t)$ denote the queue backlog of node $i$ corresponding to destination $d$ at time slot $t$,
i.e. $Q_i^d(t)$ denotes the number of commodity $d$ packets in node $i$ at time slot $t$. 
Any data that is successfully delivered to its destination will exit the network and hence, $Q_d^d(t) = 0$ for all $d \in \mathcal{D}$ and all time slots $t$.
We define $\boldsymbol{Q}(t)=[Q_i^d(t)]_{d \in \mathcal{D}, i \in \Omega}$ to be the matrix (of size $|\mathcal{D}| \times N$) of all queue backlogs and $\boldsymbol{Q}^d(t)$ to be the row of this matrix corresponding to commodity $d$ packets, i.e. $\boldsymbol{Q}^d(t)=[Q_1^d(t), Q_2^d(t), \ldots, Q_N^d(t)]$.


%
%

We define a \emph{routing decision} $\mu_{ij}^d(t)$ to be the (potential) number of commodity $d$ packets whose relaying responsibility is shifted from node $i$ to node $j$ during time slot $t$. We assume each node 
transmits at most one packet during a single time slot which can be selected from any of the 
$|\mathcal{D}|$ buffers maintained at that node.  Note that $\mu_{ij}^d(t)$ forms the departure process from 
node $i$, while it is an element of the endogenous arrival to node $j$. Hence, 
\begin{align}
\label{routingdecision}
\mu_{ij}^d(t) \in \{0,1 \} \ , \ \sum_{d \in \mathcal{D}} \mu_{ij}^d(t) \le 1\ ,\ \ 
\sum_{j=0}^{N} \mu_{ij}^d(t) \le 1.
\end{align}

Here we assume a simple on-off channel model and 
we assume a perfect channel state information at every transmitter. More specifically, 
let $S_i(t)$ represent the (random) set of nodes the channel to whom from 
node $i$ at time slot $t$ is in good state. 
We refer to $S_i(t)$ as the set of \emph{potential forwarders} for node $i$ and
we assume that node $i$ has perfect knowledge of $S_i(t)$. Due to a perfect recall at any node $i$, we assume $i \in S_i(t)$ for all time $t$. We characterize the behavior of the wireless channel using the  probabilistic model of  \emph{local broadcast model} \cite{Lott00}. The local broadcast model is defined using a marginal  probability mass functions $P(S|i):= \mbox{Prob}(\{S_i(t)= S  \ \mbox{when $i$ transmits a packet at time } t \})$, $S \subseteq \Omega$, $i \in \Omega$. 
Note that, by definition, for all $S \neq S'$, successful reception at $S$ and $S'$ are mutually exclusive and $\sum_{S \subseteq \Omega} P(S|i) =1$.
We say node $i$ \emph{reaches} node $j$ (we write $i \to j$), if there exists a set of nodes $S \subseteq \Omega$ such that $j \in S$ and $P(S|i) > 0$. 

Under the simple on-off channel model considered here, the routing of packets can only occur over links in an
on state. Often there might also be
some constraints on the simultaneous activation of  the links, i.e. certain links cannot provide service at the same time. Let an activation set be a set of links which can be activated in the same slot.
We assume that at any time $t$ the collective routing decisions 
$\{\mu_{ij}^d(t)\}_{i,j \in \Omega, d \in \mathcal{D}}$ must be such that the set of links $(i,j)$ for which $\mu_{ij}^d(t)=1$, $d\in \mathcal{D}$, belong to an activation set. Letting $\Gamma$ denote the set of all such 
allowable routing decisions, the above constraints can be written as
\begin{align}\label{routingdecisionconstrained}
\mu_{ij}^d(t) \le  \mathbf{1}_{\{ j \in S_i(t) \}} \ ,   \{ \mu_{ij}^d(t) \}_{i,j \in \Omega, d \in \mathcal{D}} \in \Gamma.
\end{align}

The selection of routing decisions together with the exogenous arrivals impact the queue backlog of node~$i$ corresponding to commodity $d$ in the following manner:
\begin{eqnarray}
\label{QDynamic}
Q_i^d(t+1) = [ Q_i^d(t) - \sum_{j \in \Omega} \mu_{ij}^d (t) ] ^{+} + \sum_{j \in \Omega} \mu_{ji}^d(t) \mathbf{1}_{\{Q_j^d(t) \ge \mu_{ji}^d(t)\}} + A_i^d(t).
\end{eqnarray}

\begin{definition}
A \emph{routing policy} is a collection of routing decisions $\cup_{i,j \in \Omega} \cup_{d \in \mathcal{D}} \cup_{t=0}^{\infty} \{ \mu_{ij}^d(t) \}$ where for all $i,j \in \Omega$, $d \in \mathcal{D}$, and $\theta \in \{0,1\}$, the decisions $\{ \mu_{ij}^d(t) = \theta \}$ belong to the $\sigma$-field generated by 
$\mathcal{H}(t)=\cup_{i,j \in \Omega} \cup_{d \in \mathcal{D}} \{ Q_i^d(0), S_i(0),  \mu_{ij}^d(0), \ldots, \\ Q_i^d(t-1), S_i(t-1), \mu_{ij}^d(t-1), Q_i^d(t), S_i(t) \}$. 
\end{definition}

\begin{definition}
A routing policy $\Pi$ is said to \emph{stabilize} the network if the time average queue backlog of each node remains finite when packets are routed according to $\Pi$. 
The \emph{stability region} of the network is the set of all arrival rate matrices $\boldsymbol{\lambda}$ for which there exists a routing policy that stabilizes the network. 
\end{definition}
\begin{definition} 
A routing policy is said to be \emph{throughput optimal} if it stabilizes the network for all arrival rate matrices that belong to the interior of the stability region. 
\end{definition}

\begin{fact}[Corollary 1 in \cite{Neely09}]
\label{Neelyeps}
Let $\mathfrak{S}$ denote the stability region of the network. An arrival rate matrix $\boldsymbol{\lambda}$ is in the stability region $\mathfrak{S}$ if and only if there exists a stationary randomized routing policy that makes routing decisions $\{\tilde{\mu}_{ij}^d(t) \}_{i,j \in \Omega, d \in \mathcal{D}}$ solely based on the collection of potential forwarders at time $t$, $\{ S_i(t) \}_{i \in \Omega}$, and for which 
$$\mathbb{E} \left[ \sum_{j \in \Omega} \tilde{\mu}_{kj}^d(t) - \sum_{i \in \Omega} \tilde{\mu}_{ik}^d(t) \right] \ge \lambda_k^d, \hspace{0.1 in} \forall k \in \Omega, \forall d \in \mathcal{D}, k \neq d.$$
\end{fact}


Fact~\ref{Neelyeps} provides a linear program whose solutions always stabilize the network, but requires a full knowledge of the arrivals statistics.
In this paper, we are interested in a class of routing policies which are throughput optimal but do not require knowledge of the arrival rates.

Now we are ready to provide the main analytical results of the paper. 
For simplicity of exposition, we first consider the single destination scenario with no activation constraints. 
The generalization of the results to the multi-destination scenario with activation constraint is provided in Section~\ref{multidst}.

Before we proceed, we introduce the following notations in the interest of simplicity:
For a set $C$ of nodes, we define $A_{C}^d(t)= \sum_{i \in C} A_i^d(t)$, 
$Q_{C}^d(t)= \sum_{i \in C} Q_i^d(t)$, $\mu_{C,in}^d(t)=\sum_{j \notin C} \sum_{k \in C} \mu_{jk}^d (t)$, and  $\mu_{C,out}^d(t)=\sum_{j \in C} \sum_{k \notin C} \mu_{jk}^d (t)$.

\section{Single Destination and Orthogonal Channel Scenario}
\label{SingleDst}
In this section, we consider the single destination with orthogonal channel scenario
and provide an overview of the results.
The analysis of the results is provided in Subsection~\ref{Thm1Proof} while the generalization of the results to the multi-destination scenario with parallel transmission constraints is discussed in Section~\ref{multidst}.

Without loss of generality, we consider node $N$ to be the destination, i.e. $\mathcal{D}=\{N\}$.
Since each node maintains only one buffer (corresponding to destination $N$), we drop the
commodity superscript $d=N$ when denoting  various
random variables such as routing decisions $\mu_{ij}(t)$, $i,j \in \Omega$. Furthermore, in this section
we assume that all channels are orthogonal and there are no activation set constraints on simultaneous 
packet transmissions, i.e. $\Gamma=\Omega \times \Omega$.

\subsection{Priority-Based Routing}


In this subsection, we introduce the class of \emph{priority-based} routing policies. To define the priority-based routing policy, we need the following definitions.

A \emph{rank ordering} $R=(C_1,C_2,\ldots,C_M)$ is an ordered list of non-empty sets $C_1,C_2,\ldots,C_M$ $(1\le M \le N)$, referred to as \emph{ranking classes}, that create a partition of $\Omega=\{ 1, 2, \ldots, N \}$, 
i.e. $\cup_{i=1}^{M} C_i = \{ 1,2,\ldots,N \}$ and $C_i \cap C_j = \emptyset$, $i \neq j$. 
We denote the set of all possible rank orderings of $\{ 1, 2, \ldots, N \}$ by $\mathcal{R}$. 
Note that when $C_i$'s are singleton, $R$ reduces to a simple permutation of the nodes $\{ 1, 2, \ldots, N \}$.
Given a rank ordering $R=(C_1,C_2,\ldots,C_M)$, we write $a \prec^R b$ to indicate that node $a \in C_i$ has a lower rank than $b \in C_j$, $i < j$. We write $a \preceq^R b$, if $a \prec^R b$ or $a,b \in C_i$ for some $i$.

\begin{definition}
A \emph{priority-based} routing policy $\Pi_{\{R(t)\}}$ is a routing policy under which node $i$, at time $t$ and among its set of potential forwarders $S_i(t)$, selects a node with the lowest rank according to $R(t)$. In other words, under $\Pi_{\{R(t)\}}$, 
$\mu_{ij}(t)=1$, only when $j \in S_i(t)$ and $j \preceq^{R(t)} k$ for all $k \in S_i(t)$. 
\end{definition}

Next we introduce a class of priority-based routing policies under which $R(t)$ is chosen as a time-invariant function of $\boldsymbol{Q}(t)$, i.e. there exists a function $\pi: \mathbb{R}^N_+ \to \mathcal{R}$ such that $R(t)=\pi(\boldsymbol{Q}(t))$. 
In Subsection~\ref{MainResults}, we proceed to establish the throughput optimality of this class of routing policies.

\subsection{$f$-policy}
\label{fpolicy}

In this section, we introduce a class of priority-based routing policies each of which is associated 
with a bivariate function~$f$, hence referred to as an $f$-policy. 
Each such policy  partitions the space of queue backlogs, $\mathbb{R}^N_+$, into $\left| \mathcal{R} \right|$ routing decision cones 
to each of which a unique rank ordering of nodes $R \in \mathcal{R}$ is assigned. 
In other words, it is possible to define the mapping  $\pi_f: \mathbb{R}^{N}_{+} \to \mathcal{R}$ such that 
at any time $t$ and for all $\boldsymbol{Q}(t)$ in the cone associated with $R$, $\pi_f(\boldsymbol{Q}(t))=R$.
The specific shape of each cone (i.e. the set of its defining hyperplanes) is dictated by the corresponding function $f$.
In order to give the precise description of $f$-policy, we need the following definitions which allow us to compare rank orderings $R$ and $R'$:
 
\begin{definition}
\label{MismatchDef}
Let $R=(C_1,C_2, \ldots,C_M)$ and $R'=(C'_1,C'_2, \ldots,C'_{M'})$. We define a {\emph{mismatch}} $m:\mathcal{R} \times \mathcal{R} \to \mathbb{N}$ as
\begin{eqnarray*}
m(R,R')= \min \left\{ i \in \mathbb{N}: C_{i} \neq C'_{i} \right\}.
\end{eqnarray*}
For two rank orderings $R$ and $R'$, $m(R,R')$ compares ranking classes of $R$ and $R'$ from low to high and determines the index of the first ranking class in which they differ. 
\end{definition}

\begin{definition}
Given two rank orderings $R$ and $R'$, we say $R'$ is a \emph{refinement} of $R$ (and $R$ is a \emph{confinement} of $R'$) if $i \prec^R j$ implies that $i \prec^{R'} j$ for any $i,j \in \Omega$.
\end{definition}

\begin{definition}
Given two rank orderings $R=(C_1,C_2,\ldots,C_M)$ and $R'=(C'_1,C'_2,\ldots,C'_{M+1})$, we say $R'$ is a \emph{one-step refinement} of $R$ (and $R$ is a \emph{one-step confinement} of $R'$) with regard to ranking class $C_i$  $(1\le i\le M)$ if 

\[ \left\{ \begin{array}{ll}
      C_k = C'_k & \mbox{if  $1 \le k \le i-1$}\\
      C_i = C'_i \cup C'_{i+1} & \mbox{} \\
      C_k = C'_{k+1} & \mbox{if $i+1 \le k \le M$}                   
         \end{array} \right. . \]
The union of the sets of all one-step refinements and one-step confinements of $R$, denoted by 
$\mathcal{B}_1(R)$ and $\mathcal{B}_2(R)$ respectively,  is referred to as \emph{adjacency} of $R$ and is denoted by $\mathcal{A} (R)$.
\end{definition}

\begin{definition}
Given a bivariate function $f$, a \emph{penalty} function $\Lambda_f$ is defined on backlog vector $\boldsymbol{Q} \in \mathbb{R}^N_+$, rank ordering $R=(C_1,C_2,\ldots,C_M) \in \mathcal{R}$, and natural number $n$, $n \le M$:
\begin{eqnarray*}
\Lambda_{f}(\boldsymbol{Q},R,n)= \sum_{i=1}^{n} f(|C^{i-1}|,|C_{i}|) Q_{C_i},
\end{eqnarray*}
where $|C_0|=0$.
\end{definition}

\begin{definition}
\label{PenalizeDef}
Consider two rank orderings $R$ and $R'$ and a bivariate function $f$. We say \emph{$R$ penalizes $\boldsymbol{Q}$ less than $R'$} 
and write $R <_{\boldsymbol{Q}} R'$ if 
\begin{itemize}
	\item {$\Lambda_{f}(\boldsymbol{Q},R,m(R,R')) < \Lambda_{f}(\boldsymbol{Q},R',m(R,R'))$,}

or if 
	\item {$\Lambda_{f}(\boldsymbol{Q},R,m(R,R')) = \Lambda_{f}(\boldsymbol{Q},R',m(R,R')) \ \text{and $R$ is a one-step refinement of $R'$.}$}
\end{itemize}
\end{definition}

Let $D_{f}(R)$, $R \in \mathcal{R}_{}$, be a subset of $\mathbb{R}^{N}_{+}$ such that for all $\boldsymbol{Q} \in D_f(R)$ and all $R' \in \mathcal{A}(R)$, $R <_{\boldsymbol{Q}} R'$, i.e.
\begin{eqnarray}
\label{GenConeDef}
D_{f}(R)= \left\{ \boldsymbol{Q} \in \mathbb{R}^{N}_{+}: R <_{\boldsymbol{Q}} R' \ \text{for all} \  R' \in \mathcal{A}_{}(R) \right\} \ .
\end{eqnarray}

\begin{remarks}
Let $R$ and $R'$ be two  rank orderings and let $\gamma \in \mathbb{R}_{+}$ be a constant. If $R <_{\boldsymbol{Q}} R'$  then $R <_{\gamma \boldsymbol{Q}} R'$. In other words, $D_f(R)$ is a cone in $\mathbb{R}^{N}_{+}$.
\end{remarks}
\begin{remarks}
Due to the linearity of  $\Lambda_{f}(\cdot,R,n)$ and finiteness of $\mathcal{A}(R)$, 
the boundaries of the cone corresponding to rank
ordering $R$ consists of finitely many hyperplanes of the form
\[
\Lambda_{f}(\boldsymbol{Q},R,m(R,R')) = \Lambda_{f}(\boldsymbol{Q},R',m(R,R')), 
\] 
where $R'\in \mathcal{A}(R)$. 
\end{remarks}

\begin{lemma}
\label{existuniq}
Let bivariate function $f$ satisfy the following two conditions
\begin{itemize}
	\item {{\bf{(C1)}} For all $m \ge 0$ and $n_1,n_2 > 0$
\begin{eqnarray*}
\frac{1}{f(m,n_1+n_2)} = \frac{1}{f(m,n_1)} + \frac{1}{f(m+n_1,n_2)} \ ,
\end{eqnarray*}
}
\item {{\bf{(C2)}} For all $m \ge 0$ and $n_1,n_2 > 0$
\begin{eqnarray*}
f(m,n_1) \ge f(m+n_1,n_2).
\end{eqnarray*}
}
	\end{itemize}
Then for any $\boldsymbol{Q} \in \mathbb{R}^{N}_{+}$, there exists a unique $R \in \mathcal{R}_{}$ such that $\boldsymbol{Q} \in D_f(R)$.
\end{lemma}

\begin{IEEEproof}
The proof is given in Appendix \ref{App:existuniq}.
\end{IEEEproof}

\begin{remarks}
By Lemma~\ref{existuniq}, $\{D_f(R)\}_{R \in \mathcal{R}}$ forms a  partition of 
$\mathbb{R}^{N}_{+}$. Hence, it is meaningful to define a function $\pi_f: \mathbb{R}^{N}_{+} \to \mathcal{R}$ such that $\pi_f(\boldsymbol{Q})=R \Leftrightarrow \boldsymbol{Q} \in D_f(R)$.
\end{remarks}

Now we are ready to provide the precise definition of $f$-policy as discussed earlier.

\begin{definition}
$f$-policy is a priority-based routing policy $\Pi_{\{ R(t) \}}$ where $R(t) = \pi_f(\boldsymbol{Q}(t))$.
\end{definition}

\begin{example}
\label{Net3nodes}
Consider a network of three nodes as given in Fig. \ref{fig:ExampleNet3}.
Let $\mathcal{R}$ be the set of all rank orderings of $\{ 1,2,3 \}$, and $f(m,n)=\frac{1}{3^m(3^n-1)}$ (it is easy to show that function $f$ satisfies (C1) and (C2)).
Since node $3$ is the destination and $Q_3(t)=0$ for all time slots $t$, the space of queue backlogs can be reduced to $\mathbb{R}^2_+$.
Furthermore, it suffices to restrict $\mathcal{R}$ to the set of all rank orderings in which the first ranking class only consists of node $3$, i.e. $C_1=\{3\}$. 
Fig. \ref{fig:LinesGen} shows the structure of the cones $\{D_f(R)\}_{R \in \mathcal{R}}$.

\begin{figure}[htp]
\centering
\subfigure
[A network of three nodes]
{
		\includegraphics[width=0.4\textwidth]{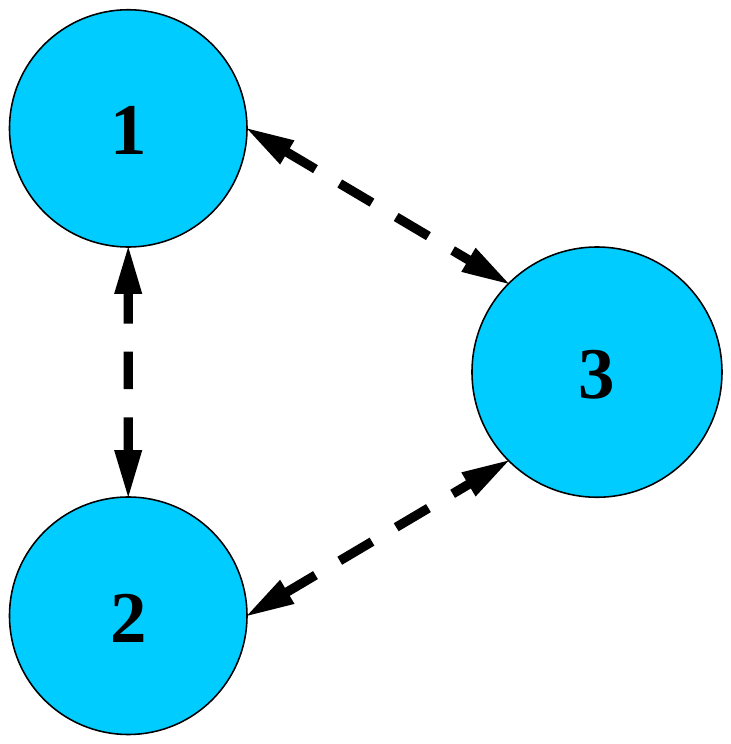}    
    \label{fig:ExampleNet3}
     }     
\subfigure
[ Structure of the cones ]
{    \includegraphics[width=0.4\textwidth]{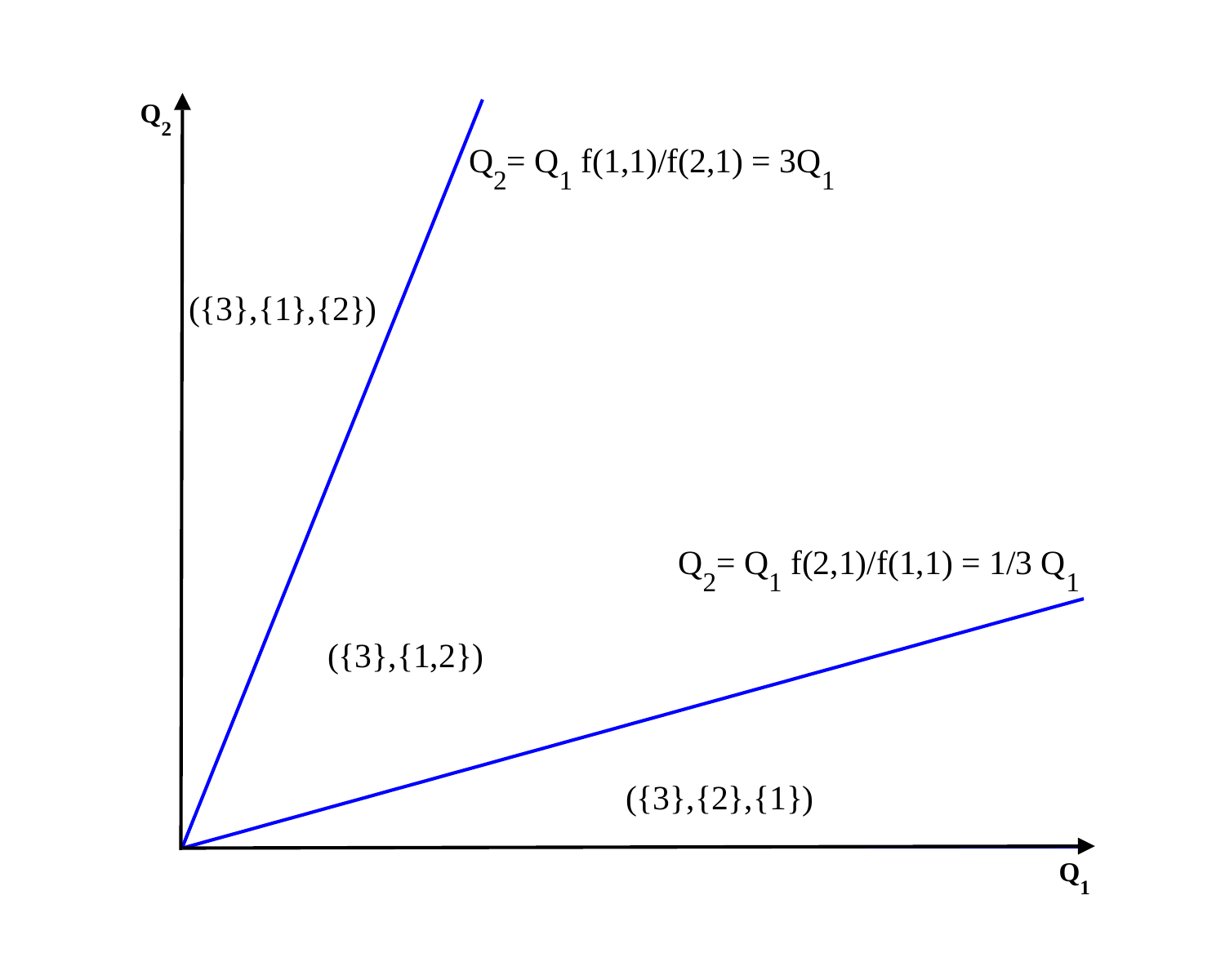}    
    \label{fig:LinesGen} 
}  
\caption{Structure of the cones for a network of three nodes }
\end{figure}
\end{example}

\begin{example}
\label{Net4nodesGen}
Consider a network of four nodes as given in Fig. \ref{fig:ExampleNet4}.
Let $\mathcal{R}$ be the set of all rank orderings of $\{ 1,2,3,4 \}$, and $f(m,n)=\frac{1}{3^m(3^n-1)}$.
Similar to Example~\ref{Net3nodes} and since $Q_4(t)=0$ for all time slots $t$, the space of queue backlogs can be reduced to $\mathbb{R}^3_+$.
Furthermore, it suffices to restrict $\mathcal{R}$ to the set of all rank orderings in which the first ranking class only consists of node $4$, i.e. $C_1=\{4\}$.  
Fig. \ref{fig:ConesGen} shows the structure of the cones $\{D_f(R)\}_{R \in \mathcal{R}}$.

\begin{figure}[htp]
\centering
\subfigure
[A network of four nodes]
{
		\includegraphics[width=0.35\textwidth]{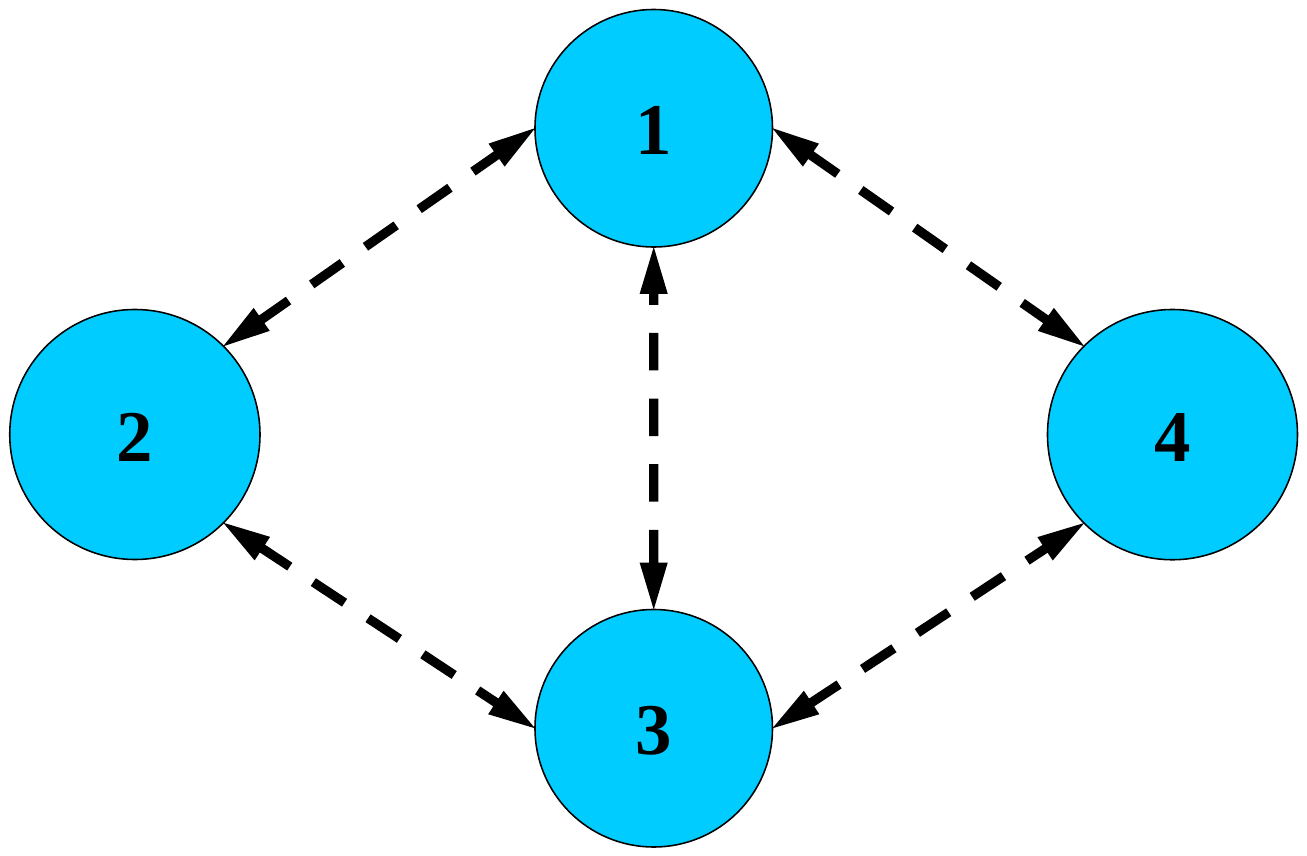}    
    \label{fig:ExampleNet4}
     }     
\subfigure
[ Structure of the cones ]
{    \includegraphics[width=0.6\textwidth]{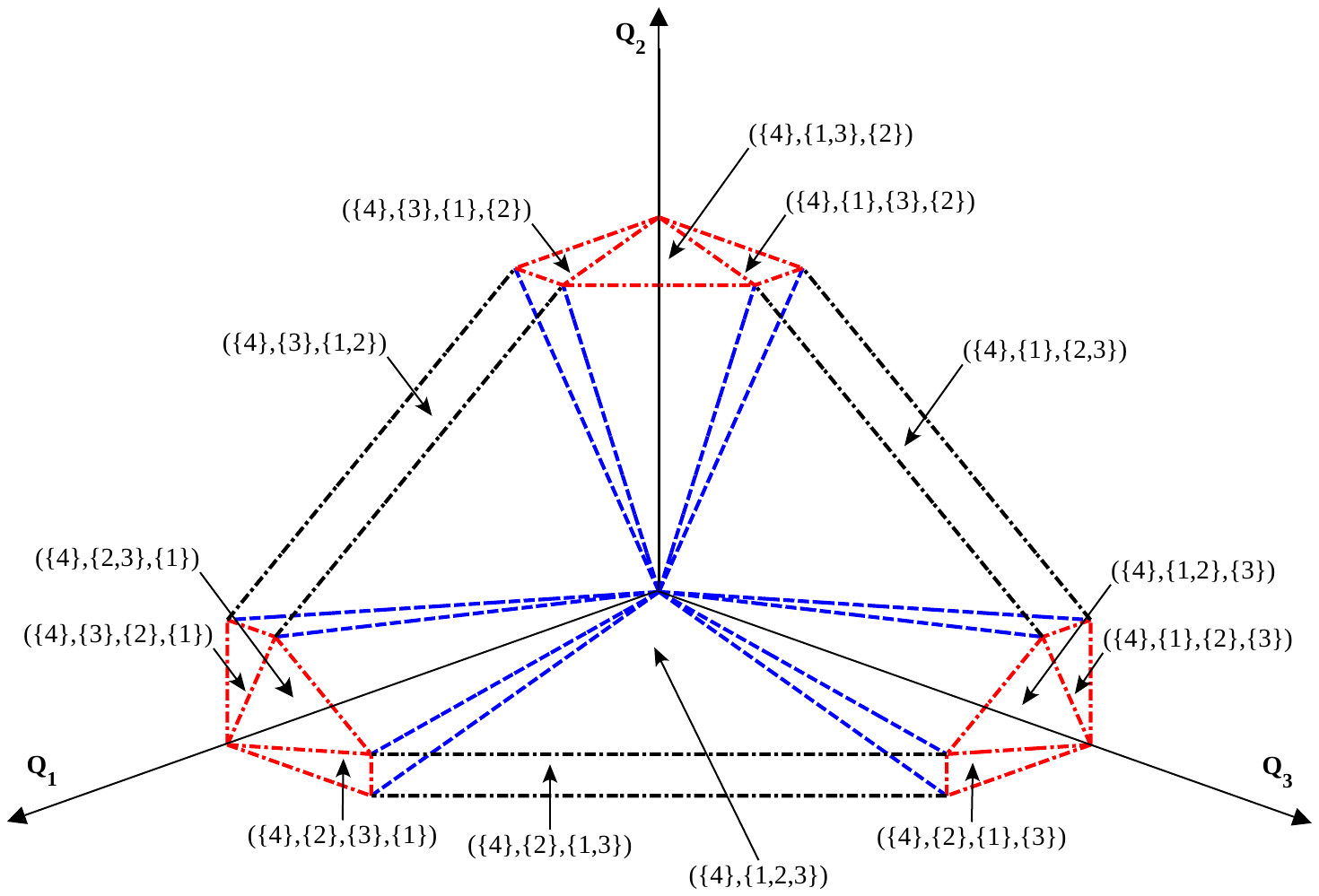}   
    \label{fig:ConesGen} 
}  
\caption{Structure of the cones for a network of four nodes }
\end{figure}

\end{example}

By construction, $f$-policy orders the nodes based only on their queue backlogs using a bivariate function $f$  independently of the topological characteristic of the network. 
In certain cases, this may cause packets to be routed away from the destination. 
In the next section, we introduce a modified version of $f$-policy, referred to as \emph{path-connected $f$-policy}, which does not allow packets to be routed away from the destination. 
The main idea behind path-connected $f$-policy is that the rank orderings are limited to those under which there exists a path from any node $i$ to the destination through the nodes with lower or the same rank as $i$. 
The precise description of path-connected $f$-policy is provided in the next section.

\subsection{Path-connected $f$-policy}
In order to give a detailed description of path-connected $f$-policy, we have to define a path-connected rank ordering.

\begin{definition}
A rank ordering $R$ is referred to as \emph{path-connected} if for each node $i$ there exist distinct nodes  $j_1,j_2,\ldots,j_l$ such that $i \to j_1 \to j_2 \to \ldots \to j_l \to N$ and $j_n \preceq^{R} i$ for all $1 \le n \le l$.
\end{definition}

The set of all path-connected rank orderings is denoted by $\mathcal{R}^{c}$, $\mathcal{R}^c \subseteq \mathcal{R}$.
Let $\mathcal{A}^{c} (R)\subseteq \mathcal{A} (R)$ be the union of the sets of all path-connected one-step refinements and one-step confinements of $R$, denoted by $\mathcal{B}_1^c(R)$ and $\mathcal{B}_2^c(R)$ respectively.
We define $D_f^c(R)$, $R \in \mathcal{R}^c$, as
\begin{eqnarray*}
D^{c}_{f}(R)= \left\{ \boldsymbol{Q} \in \mathbb{R}^{N}_{+}: R <_{\boldsymbol{Q}} R' \ \text{for all} \  R' \in \mathcal{A}^{c}(R) \right\}.
\end{eqnarray*}

\begin{definition}
The network is said to be \emph{connected} if for each node $i$ there exist nodes  $i_1,i_2,\ldots,i_l$ such that $i \to i_1 \to i_2 \to \ldots \to i_l \to N$. 
\end{definition}

Next lemma renders the set of cones as a partition of $\mathbb{R}^{N}_{+}$. 
\begin{lemma}
\label{existuniq_critical}
Assume the network is connected.{\footnote{If a node has no path to the destination, it cannot sustain any traffic and can be ignored without loss of generality.}}
If bivariate function $f$ satisfies conditions (C1) and (C2), then for all $\boldsymbol{Q} \in \mathbb{R}^{N}_{+}$, there exists a unique $R \in \mathcal{R}^{c}$ such that $\boldsymbol{Q} \in D^c_f(R)$.
\end{lemma}

\begin{IEEEproof}
The proof is given in Appendix \ref{App:existuniq}.
\end{IEEEproof}

In other words, $\{D^c_f(R)\}_{R \in \mathcal{R}^{c}}$ is the set of cones that partition $\mathbb{R}^{N}_{+}$ and it is possible to define a function $\pi^c_f: \mathbb{R}^{N}_{+} \to \mathcal{R}^c$ such that $\pi^c_f(\boldsymbol{Q})=R  \Leftrightarrow \boldsymbol{Q} \in D^c_f(R)$. 

\begin{definition}
A priority-based routing policy $\Pi_{\{ R(t) \}}$ is said to be a \emph{path-connected $f$-policy}
if $R(t) = \pi^c_f(\boldsymbol{Q}(t))$.
\end{definition}

\begin{example}
\label{Net4nodesCrit}
Consider the network of four nodes given in Example \ref{Net4nodesGen}.
Note that $(\{4\},\{2\},\{1\},\{3\})$, $(\{4\},\{2\},\{3\},\{1\})$, and $(\{4\},\{2\},\{1,3\})$, are not path-connected. 
Figure \ref{fig:ConesCrit} shows the structure of the cones $\{D^c_f(R)\}_{R \in \mathcal{R}^c}$ where $\mathcal{R}^{c}$ is the set of all path-connected rank orderings of $\{ 1,2,3,4 \}$ in which $C_1=\{4\}$ and $f(m,n)=\frac{1}{3^m(3^n-1)}$. Note the difference with Fig. \ref{fig:ConesGen} depicting $\{D_f(R)\}_{R \in \mathcal{R}}$.

\begin{figure}[htp]
\centering
    \includegraphics[width=0.6\textwidth]{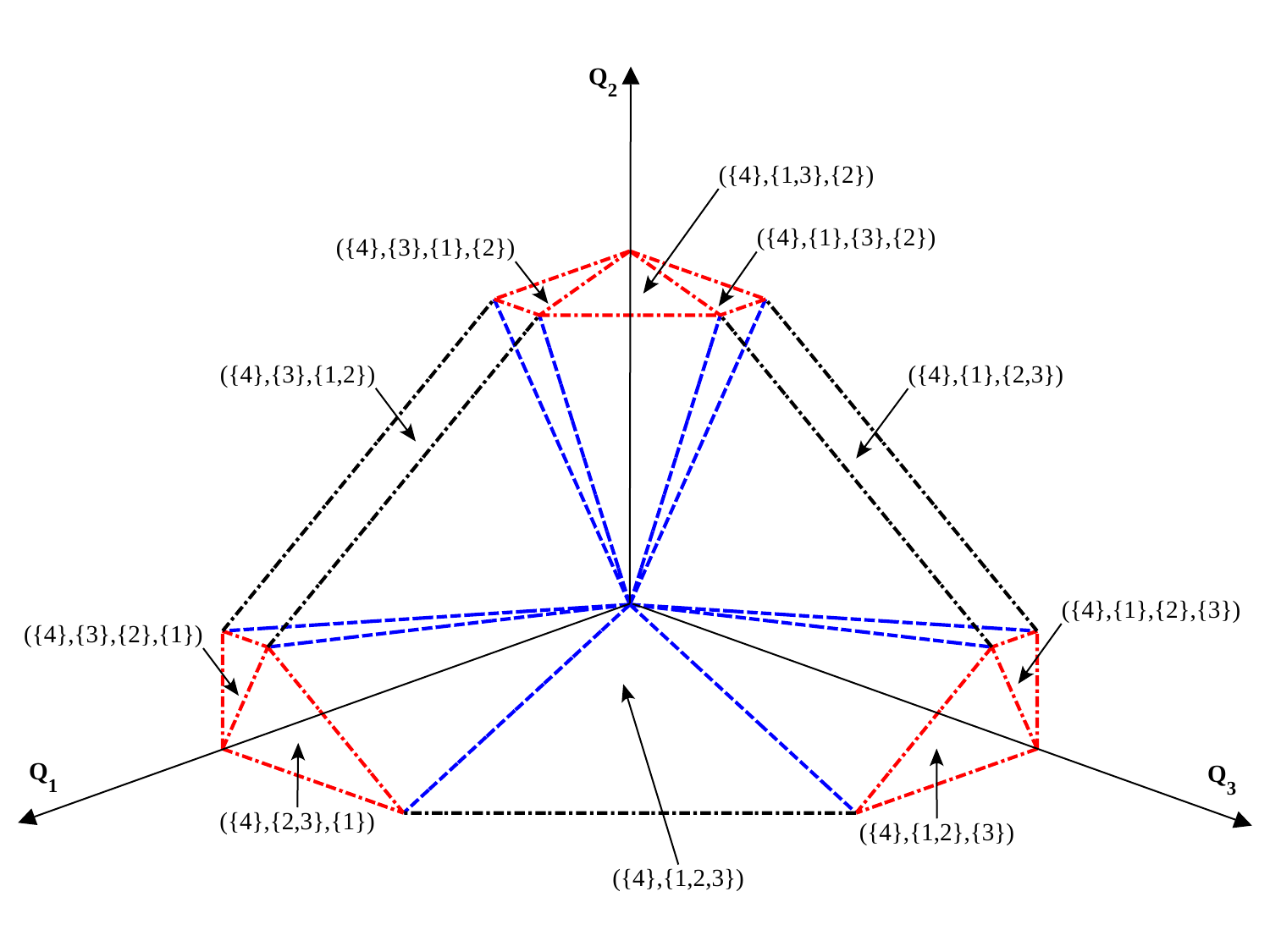}    
    \label{fig:ConesCrit} 
  
\caption{Structure of the path-connected cones for the network of Example \ref{Net4nodesGen} }
\end{figure}

\end{example}

Next we state the main results of this paper.

\subsection{Overview of the Results} 
\label{MainResults}
\begin{theorem}
\label{fpolicyopt}
Let $f$ be a bivariate function that satisfies conditions (C1) and (C2). Then the associated $f$-policy (path-connected $f$-policy) is throughput optimal.
\end{theorem}

Theorem~\ref{fpolicyopt} introduces a new class of throughput optimal routing policies. The sketch of the proof is provided in Subsection~\ref{Thm1Proof}, with the details provided in the appendix.

\begin{definition}
Let $\Pi_{\{ R(t) \}}$ and $\Pi'_{\{R'(t)\}}$ be two priority-based routing policies. We say $\Pi'_{\{R'(t)\}}$ \emph{respects} $\Pi_{\{ R(t) \}}$ if $R'(t)$ is a refinement of $R(t)$ for all time slots $t$. 
\end{definition}

\begin{theorem}
\label{respectopt}
Suppose $\Pi_{\{ R(t) \}}$ is a priority-based routing policy that is throughput optimal. Any priority-based routing policy that respects $\Pi_{\{ R(t) \}}$ is also throughput optimal. 
\end{theorem}

Note that Theorem~\ref{respectopt} enables the proof of throughput optimality of specific routing policies. For example, in Subsection~\ref{backORCD}, Theorems~\ref{fpolicyopt} and \ref{respectopt} are used to prove the throughput optimality of two known routing policies, backpressure~\cite{TassEph92} and ORCD~\cite{parul07}. The proof of Theorem~\ref{respectopt} is fairly straight forward and is given in Appendix \ref{App:respectopt}.



\subsection{Throughput Optimality of $f$-policy}
\label{Thm1Proof}

In this section, we assume that routing decisions $\{ \mu^*_{ij}(t) \}_{i,j \in \Omega}$, are made under an $f$-policy for which $f$ is a bivariate function satisfying conditions (C1) and (C2). 
In this setting, we prove that $f$-policy is throughput optimal. The proof is based on the following corollary to \emph{Foster-Lyapunov Theorem}.    

\begin{fact} [Lemma 4.1 in \cite{NeelyTass06}]
\label{LyapStable} 
Let $L^{*}: \mathbb{R}^N_{+} \to \mathbb{R}_+$ be a Lyapunov function. If there exist constants $B>0$, $\epsilon>0$, such that for all time slots $t$ we have:
$$\mathbb{E} \left[ L^{*}_{}(\boldsymbol{Q}(t+1)) - L^{*}_{}(\boldsymbol{Q}(t)) | \boldsymbol{Q}(t) \right] \le B - \epsilon \sum_{k=1}^{N} Q_k(t),$$
then the network is stable, i.e. the time average queue backlog of each node remains finite.
\end{fact}

To prove Theorem~\ref{fpolicyopt}, we identify a class of Lyapunov functions that under the corresponding $f$-policy satisfy the conditions of Fact \ref{LyapStable} for all arrival rate vectors $\boldsymbol{\lambda} \in int(\mathfrak{S})$.
In particular, we construct a piece-wise Lyapunov function, $L^{*}_{f}: \mathbb{R}^{N}_{+} \to \mathbb{R}_{+}$, by assigning to each cone $D_f(R)$, $R=(C_1,C_2,\ldots,C_M)$, a quadratic function of the queue backlogs:
\begin{eqnarray*}
L_{f}(\boldsymbol{Q},R)= \sum_{i=1}^M f(|C^{i-1}|,|C_{i}|) Q^2_{C_i}.
\end{eqnarray*}
Since the collection of cones form a partition of $\mathbb{R}^N_+$, we can combine the above quadratic functions to arrive at a piece-wise quadratic function
\begin{eqnarray}
\label{Lyapunov}
L^{*}_{f}(\boldsymbol{Q})=L_{f}(\boldsymbol{Q},\pi_f(\boldsymbol{Q}))= \sum_{R \in \mathcal{R}_{}} L_{f}(\boldsymbol{Q},R) \mathbf{1}_{\{\boldsymbol{Q} \in D_f(R)\}}.
\end{eqnarray}

\begin{lemma}
\label{contdiff}
$L^*_f(\cdot)$ is continuous and differentiable. 
\end{lemma}

Note that the continuity and differentiability of $L^*_f(\cdot)$ follow 1) the continuity and differentiability of the construction of $L_{f}(\boldsymbol{\cdot},R)$ inside the cone corresponding to $R$, as well 
as 2) the construction of penalty function on the separating hyperplanes at the boundary 
of $D_f(R)$. The details are given in Appendix \ref{App:contdiff}.



Next we provide the main steps in showing $L^*_f$ has a negative expected drift.
Let us consider the Lyapunov drift when $\boldsymbol{Q}(t) \in D_f(R)$ for some $R=(C_1,C_2, \ldots, C_M) \in \mathcal{R}$. 
By Lemma~\ref{contdiff}, $L^*_f(\cdot)$ is continuous and differentiable. Thus, we can write $L^*_f(\boldsymbol{Q}(t+1))$ in terms of its first-order Taylor expansion around $L^*_f(\boldsymbol{Q}(t))$ and we obtain 
{\allowdisplaybreaks
\begin{eqnarray}
\label{LyapDrift_01}
\nonumber
\lefteqn{L^{*}_{f}(\boldsymbol{Q}(t+1)) - L^{*}_{f}(\boldsymbol{Q}(t))} \\
\nonumber
&=& (\boldsymbol{Q}(t+1) - \boldsymbol{Q}(t)) \cdot \nabla L^{*}_{f}(\boldsymbol{Q}(t)) + O(\left\| \boldsymbol{Q}(t+1) - \boldsymbol{Q}(t) \right\|^2) \\
\nonumber
&=& \sum_{i=1}^M f(|C^{i-1}|,|C_{i}|) 2 Q_{C_i}(t) (Q_{C_i}(t+1) - Q_{C_i}(t)) + O(\left\| \boldsymbol{Q}(t+1) - \boldsymbol{Q}(t) \right\|^2) \\
\nonumber
&=& \sum_{i=1}^M f(|C^{i-1}|,|C_{i}|) \left[ Q^2_{C_i}(t+1) - Q^2_{C_i}(t) - (Q_{C_i}(t+1) - Q_{C_i}(t))^2 \right] + O(\left\| \boldsymbol{Q}(t+1) - \boldsymbol{Q}(t) \right\|^2)\\
\nonumber
&=& \sum_{i=1}^M f(|C^{i-1}|,|C_{i}|) \left[ Q^2_{C_i}(t+1) - Q^2_{C_i}(t) \right] + O(\left\| \boldsymbol{Q}(t+1) - \boldsymbol{Q}(t) \right\|^2)\\
\nonumber
&\stackrel{(a)}{\le}& B_f - 2 \sum_{i=1}^M f(|C^{i-1}|,|C_{i}|) Q_{C_i}(t) \left( \mu^*_{C_i,out}(t) - \mu^*_{C_i,in}(t) - A_{C_i}(t) \right) + O(\left\| \boldsymbol{Q}(t+1) - \boldsymbol{Q}(t) \right\|^2)\\
&\stackrel{(b)}{\le}& B_f - 2 \sum_{i=1}^M f(|C^{i-1}|,|C_{i}|) Q_{C_i}(t) \left( \tilde{\mu}_{C_i,out}(t) - \tilde{\mu}_{C_i,in}(t) - A_{C_i}(t) \right) + O(\left\| \boldsymbol{Q}(t+1) - \boldsymbol{Q}(t) \right\|^2),
\end{eqnarray} }where $B_f$ is a constant bounded real number, $\{\tilde{\mu}_{ij}(t)\}_{i,j \in \Omega}$ are routing decisions made according to the stabilizing randomized rule given in Fact~\ref{Neelyeps}, and inequalities $(a)$ and $(b)$ follow respectively from 
Lemmas~\ref{queuedynamic_ineq} and \ref{routemax} below.

\begin{lemma} 
\label{queuedynamic_ineq}
Let $R=(C_1,C_2, \ldots, C_M) \in \mathcal{R}$ and $\boldsymbol{Q}(t) \in D_f(R)$. We have
\begin{eqnarray*}
Q^2_{C_i}(t+1) - Q^2_{C_i}(t) \le \beta_f - 2 Q_{C_i}(t) (\mu^*_{C_i,out}(t) - \mu^*_{C_i,in}(t) - A_{C_i}(t)),
\end{eqnarray*}
where $\beta_f$ is a constant bounded real number.
\end{lemma}

\begin{IEEEproof}
The proof is given in Appendix \ref{App:queuedynamic_ineq}.
\end{IEEEproof}


\begin{lemma}
\label{routemax} 
Let $R=(C_1, C_2, \ldots, C_M) \in \mathcal{R}$, $\boldsymbol{Q}(t) \in D_f(R)$, and let $\{\mu^*_{ij} (t)\}_{i,j \in \Omega}$ represent routing decisions made under an $f$-policy. For any collection of routing decisions $\{\mu_{ij} (t) \}_{i,j \in \Omega}$, we have
\begin{eqnarray}
\label{routing01}
\sum_{i=1}^M f(|C^{i-1}|,|C_i|) Q_{C_i}(t) (\mu^*_{C_i,out}(t) - \mu^*_{C_i,in}(t)) \ge \sum_{i=1}^M f(|C^{i-1}|,|C_i|) Q_{C_i}(t) (\mu_{C_i,out}(t) - \mu_{C_i,in}(t)).
\end{eqnarray}
\end{lemma}

\begin{IEEEproof}
The proof is given in Appendix \ref{App:queuedynamic_ineq}. 
\end{IEEEproof}

Since $\boldsymbol{\lambda} \in int(\mathfrak{S})$, there exists a positive vector $\boldsymbol{\epsilon}$ (vector of length $N$ with all elements equal to $\epsilon$, $\epsilon > 0$) such that $\boldsymbol{\lambda} + \boldsymbol{\epsilon} \in \mathfrak{S}$. Thus, from Fact~\ref{Neelyeps} 
\begin{eqnarray}
\label{epseq}
\mathbb{E} \left[ \tilde{\mu}_{C_i,out}(t) - \tilde{\mu}_{C_i,in}(t) - A_{C_i}(t) | \boldsymbol{Q}(t) \right] \ge \epsilon.
\end{eqnarray}
Now taking expectation from both sides of (\ref{LyapDrift_01}) and using (\ref{epseq}) we obtain,
\begin{eqnarray}
\label{proofThm1_01}
\mathbb{E} \left[ L^{*}_{f}(\boldsymbol{Q}(t+1)) - L^{*}_{f}(\boldsymbol{Q}(t)) | \boldsymbol{Q}(t) \right]
&\le& B_f - 2 \epsilon \sum_{i=1}^M f(|C^{i-1}|,|C_{i}|) Q_{C_i}(t) + O(\left\| \boldsymbol{Q}(t+1) - \boldsymbol{Q}(t) \right\|^2). 
\end{eqnarray}
Since $\left\| \boldsymbol{Q}(t+1) - \boldsymbol{Q}(t) \right\|$ is bounded, there exists a constant, say $B'_f$, such that 
$B_f + O(\left\| \boldsymbol{Q}(t+1) - \boldsymbol{Q}(t) \right\|^2) \le B'_f$ for all time slots $t$.
Moreover, property (C2) of function $f$ implies that
\begin{eqnarray}
\label{proofThm1_02}
f(0,|C_{1}|) \ge f(|C^{1}|,|C_{2}|) \ge \cdots \ge f(|C^{M-1}|,|C_{M}|) \ge f(|C^{M}|,1) = f(N,1).
\end{eqnarray}
Therefore, we can rewrite (\ref{proofThm1_01}) as
\begin{eqnarray*}
\mathbb{E} \left[ L^{*}_{f}(\boldsymbol{Q}(t+1)) - L^{*}_{f}(\boldsymbol{Q}(t)) | \boldsymbol{Q}(t) \right] \le B'_f - \epsilon' \sum_{k=1}^N Q_{k}(t),
\end{eqnarray*}
where $\epsilon'=2 \epsilon f(N,1)$. 
Now from Fact \ref{LyapStable}, the proof of Theorem~\ref{fpolicyopt} is complete.   

Note that the proof of throughput optimality for path-connected $f$-policy follows similar lines above and is provided in Appendix~\ref{criticalfpolicy_opt}.


\section{Generalization: Multi-Destination Constrained $f$-policy}
\label{multidst}



In this section, we introduce \emph{multi-destination constrained $f$-policy} as a 
generalization of $f$-policy in a multi-destination scenario with parallel transmission constraints.
Next we provide a precise definition of multi-destination $f$-policy.

Suppose $\boldsymbol{Q}^d(t) \in D_f(R^d)$, $R^d=(C_1^d,C_2^d,\ldots,C_{|R^d|}^d)$, $d \in \mathcal{D}$, and let $C^{i-1,d} = \cup_{j=1}^{i-1} C_j^d$.
Multi-destination constrained $f$-policy is defined as to select routing decisions $\{\mu_{ij}^d(t)\}$ such that 
for any global channel state $\{S_i(t)\}_{i\in\Omega}$, they maximize
\[
\sum_{d \in {\mathcal{D}}} \sum_{i=1}^{|R^d|} \sum_{k \in C_i^d} \sum_{j=1}^{|R^d|} \sum_{l \in C_j^d} \mu_{kl}^d(t) \left [ f(|C^{i-1,d}|,|C_{i}^d|) Q_{C^d_i}(t) - f(|C^{j-1,d}|,|C_{j}^d|) Q_{C^d_j}(t) \right ], 
\]
while satisfying (\ref{routingdecision}) and (\ref{routingdecisionconstrained}). 
Note that due to the global nature of the activation set constraints, the policy does not have the decentralized 
structure of the $f$-policy. 

\begin{theorem}
\label{fmltdstopt}
Let $f$ be a bivariate function that satisfies conditions (C1) and (C2). Then the associated multi-destination constrained $f$-policy is throughput optimal.
\end{theorem}

\begin{IEEEproof}
The proof is very similar to the proof of Theorem~\ref{fpolicyopt} provided in Subsection~\ref{Thm1Proof}.
Similar to (\ref{Lyapunov}), we define a piece-wise quadratic function $L_f^*$ as follows:
\begin{align}
L_f^*(\boldsymbol{Q}) = \sum_{d \in \mathcal{D}} L(\boldsymbol{Q}^d,\pi_f(\boldsymbol{Q}^d)).
\end{align}
Let $\{ \mu_{ij}^{*d}(t) \}_{i,j \in \Omega, d \in \mathcal{D}}$ represent  routing decisions made under a multi-destination constrained $f$-policy.
Let us consider the Lyapunov drift when $\boldsymbol{Q}^d(t) \in D_f(R^d)$, $R^d=(C_1^d,C_2^d,\ldots,C_{|R^d|}^d)$, $d \in \mathcal{D}$.
Following similar steps as that of the proof of Theorem~\ref{fpolicyopt}, we obtain
\begin{align}
\label{mltdstdrift}
\nonumber
\lefteqn{ \mathbb{E} \left [ L^{*}_{f}(\boldsymbol{Q}(t+1)) - L^{*}_{f}(\boldsymbol{Q}(t)) | \boldsymbol{Q}(t) \right ]} \\
&\le& B_f - 2 \sum_{d \in \mathcal{D}} \sum_{i=1}^{|R^d|} f(|C^{i-1,d}|,|C_{i}^d|) Q_{C_i^d}^d(t) \mathbb{E} \left[ \mu^{*d}_{C_i^d,out}(t) - \mu^{*d}_{C_i^d,in}(t) - A_{C_i^d}^d(t) | \boldsymbol{Q}(t) \right] + o(\left\| \boldsymbol{Q}(t+1) - \boldsymbol{Q}(t) \right\|),
\end{align}
where $B_f$ is a constant bounded real number.
However, the term
\begin{align}
\nonumber
\lefteqn{ \sum_{d \in \mathcal{D}} \sum_{i=1}^{|{R}^d|} f(|C^{i-1,d}|,|C_i^d|) Q_{C_i^d}^d(t) \left( \mu_{C_i^d,out}^d(t) - \mu_{C_i^d,in}^d(t) \right) } \\
\nonumber
&=  \sum_{d \in \mathcal{D} } \sum_{i=1}^{|R^d|} \sum_{k \in C_i^d} \sum_{j=1}^{|R^d|} \sum_{l \in C_j^d} \mu_{kl}^d(t) \left[ f(|C^{i-1,d}|,|C_{i}^d|) Q_{C^d_i}(t) - f(|C^{j-1,d}|,|C_{j}^d|) Q_{C^d_j}(t) \right]
\end{align} 
is maximized by the multi-destination constrained $f$-policy for any global channel state $\{S_i(t)\}_{i\in\Omega}$. Hence, the negative drift term in (\ref{mltdstdrift}) is bounded by the negative drift under any other set of routing decisions, including the stabilizing randomized rule. 
Now from Facts~\ref{Neelyeps} and \ref{LyapStable}, the proof of Theorem~\ref{fmltdstopt} is complete.

\end{IEEEproof}

As a special case of routing in constrained queueing networks, scheduling for single-hop networks (e.g. wireless uplinks and downlinks) and input-queued switches have also been of great interest \cite{Tassiulas93,Neely08,Markakis09,McKeown99,Shah05,Ross09}.
Next, we show that throughput optimal scheduling for input-queued switches is a special case of 
our framework and discuss the $f$-scheduling. 
In particular, we specialize our result to derive a class of throughput optimal scheduling policies for input-queued switches which will be compared with some of the existing scheduling policies.

\subsection{Input-Queued Switches} 
\label{IQswitch}

\begin{figure}[htp]
\centering
\includegraphics[width=0.6\textwidth]{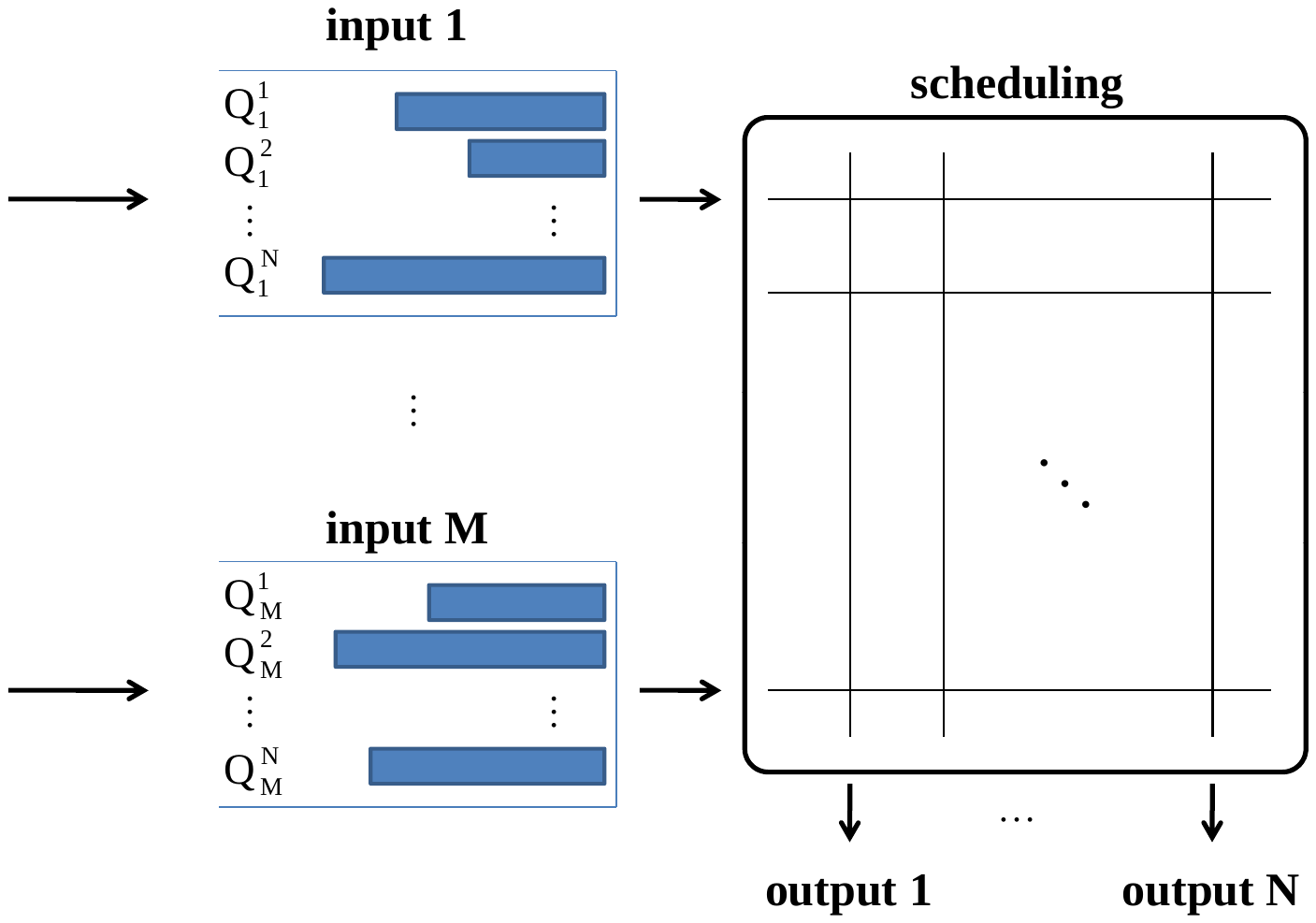}    
\caption{$M \times N$ input-queued switch.}
\label{fig:MNswitch}
\end{figure}

Consider the  \emph{input-queued switch} studied
 in \cite{McKeown99} and as depicted in Fig.~\ref{fig:MNswitch}.
 A \emph{scheduling decision} $\eta_{i}^d(t)$ is defined to be the (potential) number of packets 
sent from input $i$, $i \in \mathcal{I}$ to output $d$, $d \in \mathcal{D}$ during time slot $t$.
In a crossbar switch, each input can send to at most one output and each output can receive from at most one input and hence,  
\begin{eqnarray}
\label{mappingdecision}
\eta_{i}^d(t) \in \{0,1 \} \ ,  \ \sum_{d=1}^{N} \eta_{i}^d(t) \le 1 \ , \ \sum_{i=1}^{M} \eta_{i}^d(t) \le 1.
\end{eqnarray}


This is nothing but a single-hop example of the setup introduced in Section~\ref{ProblemFormulation} where the nodes in the network 
are partitioned into $M$ inputs labeled by $\mathcal{I}=\{1,2,\ldots,M\}$ 
to $N$ outputs (destinations) labeled by $\mathcal{D}=\{1,2,\ldots,N\}$, with $\lambda^d_{d'} = 0$ for all
$d, d' \in \mathcal{D}$ and a deterministic and fully connected local broadcast model $P(\mathcal{D} |i)=1$ for all 
$i \in \mathcal{I}$.  In this setup, choice of routing decisions $\mu_{id}^d(t)$ are equivalent to scheduling decisions $\eta_{i}^d(t)$, while the set of allowable routing decisions $\Gamma$ is the space of all permutation matrices. This means that, in the input-queued switch problem, (\ref{routingdecision}) and (\ref{routingdecisionconstrained}) reduce to (\ref{mappingdecision}). Now for any bivariate function $f$,  we introduce a class of scheduling policies each of which is constructed using our proposed framework, hence referred to as an $f$-scheduling. 
Let $\mathcal{R}$ denote the set of all possible rank orderings of $\mathcal{I}$.
$f$-scheduling  partitions the space of queue backlogs corresponding to each destination, $\mathbb{R}^{M}_+$, into $\left| \mathcal{R} \right|$ scheduling decision cones 
to each of which a unique rank ordering $R \in \mathcal{R}$ is assigned. 
Similar to (\ref{GenConeDef}) in Subsection~\ref{fpolicy}, we can define $D_f (R)$, $R \in \mathcal{R}$, such that
$\{D_f(R)\}_{R \in \mathcal{R}}$ forms a partition of $\mathbb{R}^{M}_+$. 
Suppose $R^d=(C_1^d,C_2^d,\ldots,C_{|R^d|}^d) \in \mathcal{R}$ and $\boldsymbol{Q}^d(t) \in D_f(R^d)$ for all $d \in \mathcal{D}$.
Then $f$-scheduling selects scheduling decisions $\{\eta_{i}^d(t)\}_{i \in \mathcal{I}, d \in \mathcal{D}}$ in order to maximize 
\begin{equation}
\sum_{d=1}^N \sum_{i=1}^{|R^d|} f(|C^{i-1,d}|,|C_i^d|) Q_{C_i^d}^d(t) \eta_{C_i^d}^d(t), \label{f-schedule}
\end{equation}
where $Q_{C_i^d}^d(t) = \sum_{k \in C_i^d} Q_{k}^d(t)$ and  $\eta_{C_i^d}^d(t) = \sum_{k \in C_i^d} \eta_{k}^d(t)$.

\begin{corollary}
\label{fschedulingopt}
Let $f$ be a bivariate function that satisfies conditions (C1) and (C2). Then the associated $f$-scheduling given by (\ref{f-schedule}) is throughput optimal.
\end{corollary}


\section{$f$-policy and the Design for Throughput Optimality}
\label{Design}

\subsection{Examples: The Structure of Lyapunov Function}
\label{Exmf}

\begin{figure}[!h]
	\centering
 		\includegraphics[width=0.5\textwidth]{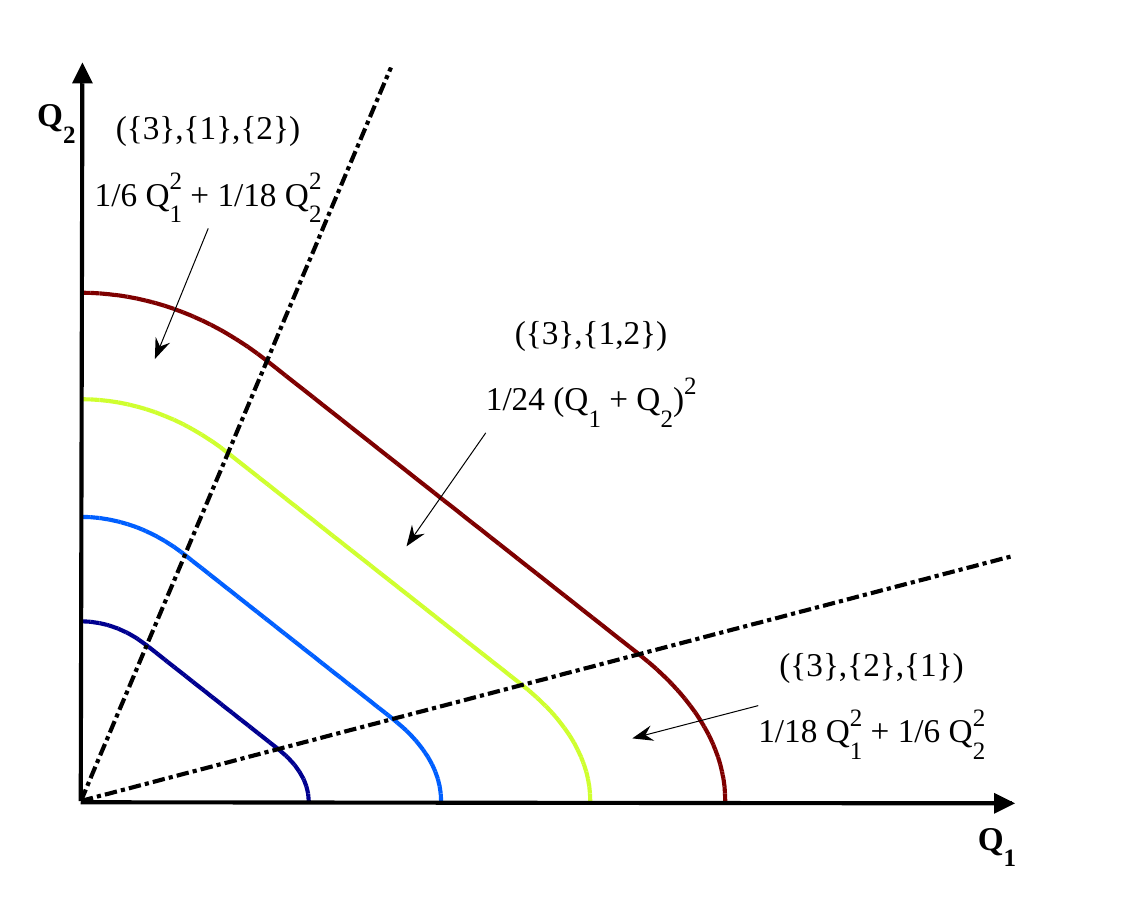}
	\caption{Structure of the Lyapunov function for the network of Example \ref{Net3nodes}}
	\label{fig:LyapContour}
\end{figure}

In Examples~\ref{Net3nodes}-\ref{Net4nodesCrit}, we considered two different networks and showed the structure of (path-connected) cones for those networks. In this section, we study the Lyapunov function as defined in (\ref{Lyapunov}) for the same networks. 
Figure~\ref{fig:LyapContour} illustrates routing decision cones $D_f(R)$ and the associated quadratic function $L(\cdot,R)$ for the network of Example~\ref{Net3nodes}.
In addition, Fig.~\ref{fig:LyapContour} shows the contour lines of the Lyapunov function.
In the central cone where nodes 1 and 2 belong to the same ranking class, the contour of the Lyapunov function is a straight line with slope 135 degrees.
The Lyapunov drift in this case is the same for all non-idling routing policies.
The contours of the Lyapunov function in the two corner cones are elliptic.
Furthermore, the contours are perpendicular to the Q axis. 
This implies that when one of the queues is close to empty, to most efficiently reduce the Lyapunov
drift, the policy tends to give a higher rank to the other queue.
In other words, the longer queue is served more often and the shorter queue is more likely
to stay away from empty state in future.
On the boundaries between the cones, i.e. lines $Q_2 = 3 Q_1$ and $Q_2 = \frac{1}{3} Q_1$, the contours are still smooth
indicating that the Lyapunov function is continuous and differentiable on the boundaries.

\begin{figure}[!h]
	\centering
 		\includegraphics[width=0.8\textwidth]{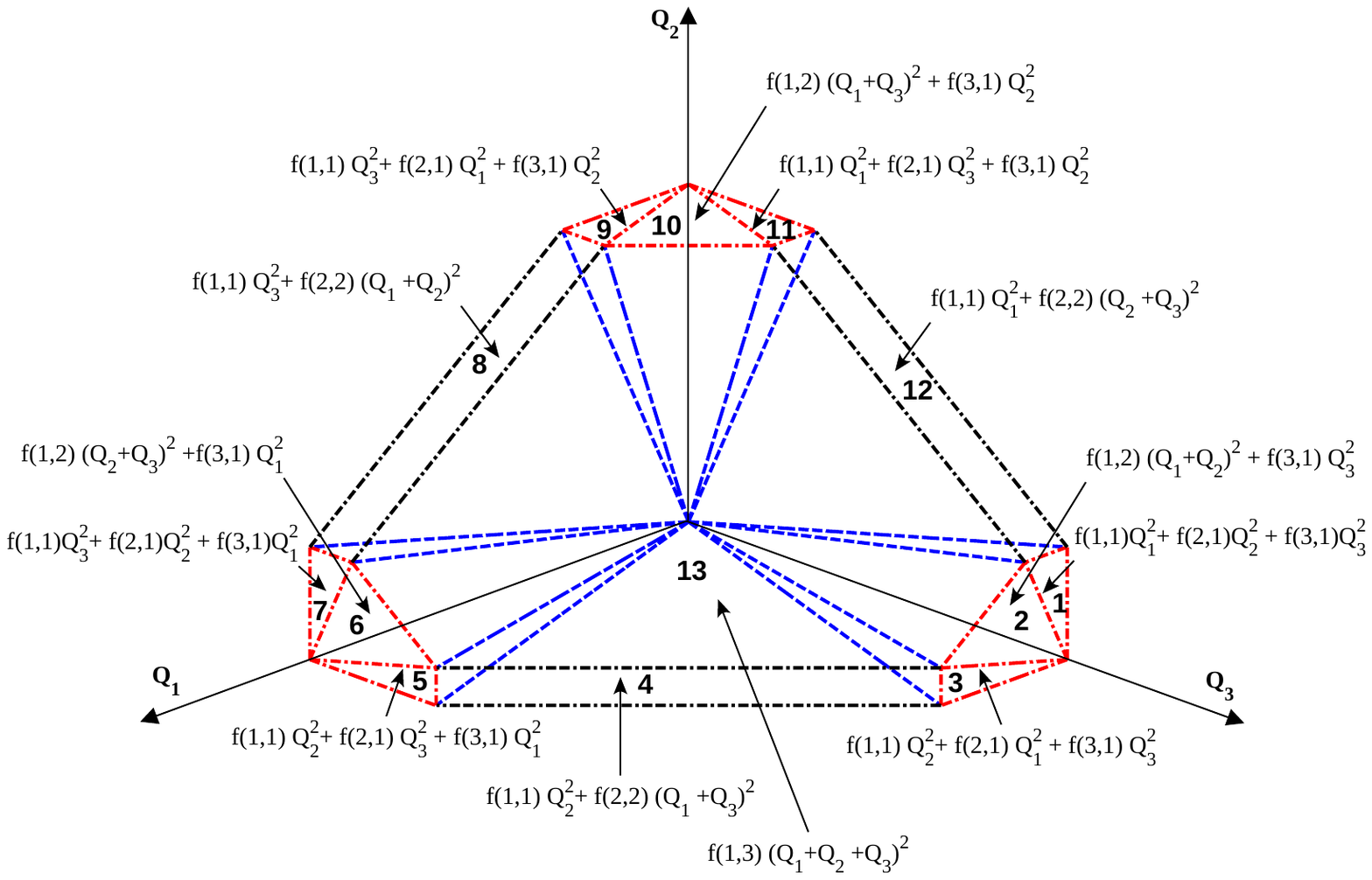}
	\caption{Structure of the Lyapunov function for the network of Example \ref{Net4nodesGen}}
	\label{fig:LyapStructure}
\end{figure}

Figure~\ref{fig:LyapStructure} illustrates routing decision cones $D_f(R)$ and the associated quadratic function $L(\cdot,R)$ for the network of Example~\ref{Net4nodesGen}.
As shown in Fig. \ref{fig:LyapStructure}, $f$-policy groups the queues based on their backlogs consistent with the cone in which the backlog state lies. 
Given such a grouping, the Lyapunov function is constructed by considering the sum of quadratic group backlogs. For instance,
in the central cone, $13$, where all nodes belong to the same ranking class, the Lyapunov function is nothing but the squared sum of all queue backlogs. It is clear that the Lyapunov drift in this case is the same (and negative from Theorem \ref{fpolicyopt}) for all non-idling routing policies. 
This property allows a routing policy to potentially deviate from backpressure decisions 
while still ensuring the throughput optimality. 
However, when one of the queues becomes relatively large in comparison to the other nodes' backlogs, the backlog vector falls in one of the cones $2$, $6$, or $10$, in which the node with large backlog is in a separate ranking class. The Lyapunov function over each of these cones is the squared queue backlog of the node with large backlog plus the squared sum of other queue backlogs. Hence, in cones $2$, $6$, or $10$, the negative expected drift is ensured only when packets are routed away from the node with disportionately large backlog. 
Similarly one can analyze the behavior of the Lyapunov function in the remaining cones. 

It is important to observe that the Lyapunov function in cones $1$, $3$, $5$, $7$, $9$, and $11$ is a weighted quadratic function and closely related to the quadratic Lyapunov function associated with backpressure routing. In contrast, the Lyapunov function in cone $13$ is closely related to the total network backlog whose negative drift can be ensured so long as the routing policy 
does not allow idling. 
In this construction, we were strongly motivated by the discussion in \cite{Meyn08} on the piece-wise structure of the \emph{value function} and its properties near the subspaces of the form $\{Q_i=0\}_{i \in I}$, for some $I \subseteq \{1,2, \ldots, N-1\}$. In particular, we have carefully constructed a smooth Lyapunov function by ``stitching" a backpressure-like behavior at the neighborhood of the subspaces $\{Q_i=0\}_{i \in I}$, i.e.\ in cones $1$, $3$, $5$, $7$, $9$, and $11$, while allowing for an arbitrary non-idling behavior far from these subspaces, e.g.\ cone $13$. 
Furthermore, function $f$ determines the size of cones associated with various rank orderings and provides flexibility in the relative occupancy time spent in non-idling cone $13$ versus backpressure-like cones $1$, $3$, $5$, $7$, $9$, and $11$. Note that, however, the flexibility in cone $13$ and non-uniqueness of the policy which ensures a negative Lyapunov drift, in effect,
prevent the proposed Lyapunov function to lead to a unique policy construction; instead, the proposed Lyapunov construction is useful either in verifying the throughput optimality of a given policy or in modification of existing policies.

\subsection{Single Commodity Multi Hop Networks: Known Routing Policies}
\label{backORCD}

Opportunistic routing for multi-hop wireless ad-hoc networks has seen recent research interest to overcome deficiencies of conventional routing \cite{Lott00,Larson01,Zorzi03,Morris05,Das05,Neely09,parul07}. 
 Under opportunistic routing, the routing decisions are made in an online manner by choosing the next relay based on the channel state realization during that time slot as well as a rank ordering of neighboring nodes. 
 Within the opportunistic routing frame work, however, the rank ordering of the neighbors, or equivalently the selection criterion when selecting the next hop, remains an important design problem. While some authors advocated for relaying packets via the neighbor with the minimum expected transmissions \cite{Larson01,Morris05}, the neighbor with minimum expected cost \cite{Lott00}, or the shortest (geographic) distance to the destination \cite{Zorzi03}, others have pointed out that when multiple streams of packets are to traverse the network, however, it might be necessary to route some packets along longer paths, especially if these paths eventually lead to links that are less congested. More precisely, \cite{Neely09}, \cite{parul07} have 1) showed that the above routing schemes fail to stabilize otherwise stabilizable traffic (see examples given in \cite{parul07}), and 2) proposed selecting the relay based on a time varying notion of congestion, such as the queue backlog differentials \cite{Neely09} or  the estimated draining time \cite{parul07}.  In other words, when multiple streams of packets are to traverse the network, it is critical for the routing solution 
 to ensure queue stability for all stabilizable traffic conditions. 
 
While one can intuitively design and propose various time variant notions of congestion, it is far from straight forward to
verify the throughput optimality of such solutions. This is mainly due to the difficulty in analyzing multi-hop, multi-server, and multi queue systems in time. In particular, the analytic guarantees usually depend on the construction of a Lyapunov function with an expected negative drift; a task far from trivial! 
 This has meant that many successful strategies are reverse engineered to be the very 
 rule for which a known Lyapunov function is ensured to have a negative expected drift. 
 In fact backpressure \cite{TassEph92} and its variants \cite{Neely09,sarkar08,Ying08,Xi06}, with quadratic 
 Lyapunov function, and randomized strategies \cite{Hassibi08} with an exponential Lyapunov function 
 remain to be the only known throughput optimal routing policies.
In this section, we use Theorems~\ref{fpolicyopt} and \ref{respectopt} to prove the throughput optimality of two known routing policies, backpressure~\cite{TassEph92} and ORCD~\cite{parul07} in a single destination with orthogonal channel scenario.

In the opportunistic variant of backpressure routing, DIVBAR \cite{Neely09}, among the set of nodes that have received a packet transmitted by node $i$, one of the nodes with the largest positive differential queue backlog is selected as the next forwarder. 
Therefore, backpressure is a priority-based routing policy $\Pi_{\{ R_{\text{b}}(t) \}}$ where $R_{\text{b}}(t)$ is a partitioning of the nodes based on their queue backlogs, with smaller backlog meaning lower rank, i.e. $Q_i(t) < Q_j(t)$ implies that $i \prec^{R_{\text{b}}(t)} j $.
This policy is shown to provide throughput optimality \cite{Neely09}. 
Here, however, we give an alternative proof which relies on Theorem~\ref{respectopt}. 
In other words, backpressure routing is proved to be throughput optimal by showing that for any bivariate function $f$ that satisfies conditions (C1) and (C2), backpressure respects $f$-policy.
More precisely, we show that if node $j$ has a lower rank than node $k$ under any $f$-policy, then $Q_j < Q_k$.
The proof is immediate using Lemma~\ref{Lemma:basic2} in Appendix \ref{Preliminary}.





In the rest of this section, we give a brief description of another congestion-based routing policy, known as ORCD \cite{parul07} and prove its throughput optimality.
In \cite{parul07}, ORCD was introduced as an alternative to backpressure routing to improve the delay performance.
However, the throughput optimality of ORCD was left as a conjecture.


ORCD is a priority-based routing policy $\Pi_{\{R_{\text{CD}}(t)\}}$ in which nodes are ordered according to a cost measure of congestion ``down the stream'' from each node $i$ denoted by $V_i(t)$. In other words, $i \prec^{R_{\text{CD}}(t)} j$ if $V_i(t) < V_j(t)$.
The congestion cost measures for nodes $i \in \Omega$ at time $t$, $V_i(t)$'s, form a vector $[V_1(t), V_2(t), \ldots , V_N(t)]$ that satisfies the following fixed point equation:
\begin{eqnarray}
\label{fixpoint0}
V_N(t) &=& 0, \\
\label{fixpoint1}
V_i(t) &=& Q_i(t) + \sum_{S \subseteq \Omega} P(S|i) \min_{j \in S} V_j(t), \ \ \ \text{for} \ i=1,2,\ldots,N-1.
\end{eqnarray}



Here, we prove the throughput optimality of ORCD by showing that ORCD respects path-connected $f$-policy corresponding to any bivariate function $f$ that satisfies condition (C1) and for all $m \ge 0$ and $n_1,n_2 > 0$ 
\begin{eqnarray}
\label{C3}
\frac{f(m,n_1)}{f(m+n_1,n_2)} \ge \frac{1}{p_{\min}},
\end{eqnarray} 
where $p_{\min}=\min \left\{ P(S|i): i \in \Omega, S \subseteq \Omega, P(S|i)>0 \right\}$. Note that, for instance, function $f(m,n)=\frac{1}{K^m(K^n-1)}$, $K \ge 1 + \frac{1}{p_{\min}}$, is such a function.
In other words, we show that ORCD respects the path-connected $f$-policy for all such $f$. Mathematically, for all $j,k \in \Omega$ such that $j \prec^{\pi^c_f(\boldsymbol{Q}(t))} k$, then $j \prec^{R_{\text{CD}}(t)} k$ as well. Let $\pi^c_f(\boldsymbol{Q}(t)) = (C_1,C_2,\ldots,C_M) \in \mathcal{R}^c$, and let $k \in C_i$ and $j \in C^{i-1}$. We consider two cases:

{\bf{Case I.}}
Node $k$ reaches a node in $C^{i-1}$. In such a case, we need to show
\begin{eqnarray}
\label{ORCDopt4}
V_k (t) \ge Q_k (t) > \frac{Q_{C^{i-1}} (t)}{p_{\min}} \ge V_j (t).
\end{eqnarray}  
%
The first inequality in (\ref{ORCDopt4}) is immediate from (\ref{fixpoint2}). The second and third inequalities follow from the arguments below.
Lemma~\ref{Lemma:basic4} in the appendix implies the second inequality in (\ref{ORCDopt4}), i.e.\
\begin{eqnarray}
\label{ORCDopt1}
Q_k (t) > \frac{f(0,|C^{i-1}|)}{f(|C^{i-1}|,1)} Q_{C^{i-1}} (t) \ge \frac{Q_{C^{i-1}} (t) }{p_{\min}}.
\end{eqnarray}
On the other hand, since $\pi^c_f(\boldsymbol{Q}(t))$ is path-connected, there exist distinct intermediate nodes $j_1, j_2, \ldots, j_l \in C^{i-1}$ such that $j \to j_1 \to j_2 \to \ldots \to j_l \to N$. Using Lemma~\ref{App:Vablemma} in the appendix recursively and noting that $V_N (t)=0$, we have the following upper bound of $V_j (t)$,
\begin{eqnarray}
\label{ORCDopt3}
V_j (t) \le \frac{Q_j (t)}{p_{\min}} + \frac{Q_{j_1} (t)}{p_{\min}} + \cdots + \frac{Q_{j_l} (t)}{p_{\min}} \le \frac{Q_{C^{i-1}} (t)}{p_{\min}},
\end{eqnarray}
which gives the last inequality in (\ref{ORCDopt4}).

{\bf{Case II.}}
Node $k$ does not reach any node in $C^{i-1}$.
Let $\hat{C}_i$ be the set of nodes in $C_i$ that reach a node in $C^{i-1}$. 
All the paths from node $k$ to the destination are through the nodes in $\hat{C}_i$, hence, $V_k(t) \ge \min_{m \in \hat{C}_i} V_m(t)$. However, from Case I, for each node $j \in C^{i-1}$ and $m \in \hat{C}_i$, $V_m(t) \ge V_j(t)$.
This completes the proof.

\subsection{Input-Queued Switches: Known Scheduling policies}

In this subsection, we compare $f$-scheduling with two known throughput optimal scheduling policies:
1) Longest Queue First (LQF) \cite{McKeown99}; and 2) Projective Cone Scheduling (PCS) \cite{Ross09}.
Before we proceed, a brief overview of LQF and PCS is provided. Let ${\boldsymbol{\eta}}(t)$ and $\boldsymbol{Q}(t)$ be respectively the matrices of the scheduling decisions and the queue backlogs, i.e. 
${\boldsymbol{\eta}}(t):=[\eta_{i}^d(t)]_{M \times N}$ and $\boldsymbol{Q}(t):=[Q_{i}^d(t)]_{M\times N}$. PCS is a 
class of scheduling policies each of which is parameterized with a positive-definite symmetric matrices (of size $M \times M$)
with negative or zero off-diagonal elements, $\boldsymbol{W}$.  Associated with matrix $\boldsymbol{W}$, PCS selects scheduling decisions $\{\eta_{i}^d(t)\}$ such that 
$\langle {\boldsymbol{\eta}}(t) , \boldsymbol{W} \boldsymbol{Q}(t)\rangle$  
is maximized. When $\boldsymbol{W}$ is set to be equal to identity, PCS coincides with LQF policy which  
selects scheduling decisions $\{\eta_{i}^d(t)\}$ such that 
$\sum_{d}\sum_i Q_i^d(t) \eta_i^d(t)$ is maximized.

\begin{figure}[htp]
\centering
\subfigure
[$f$-scheduling]
{
		\includegraphics[width=0.3\textwidth]{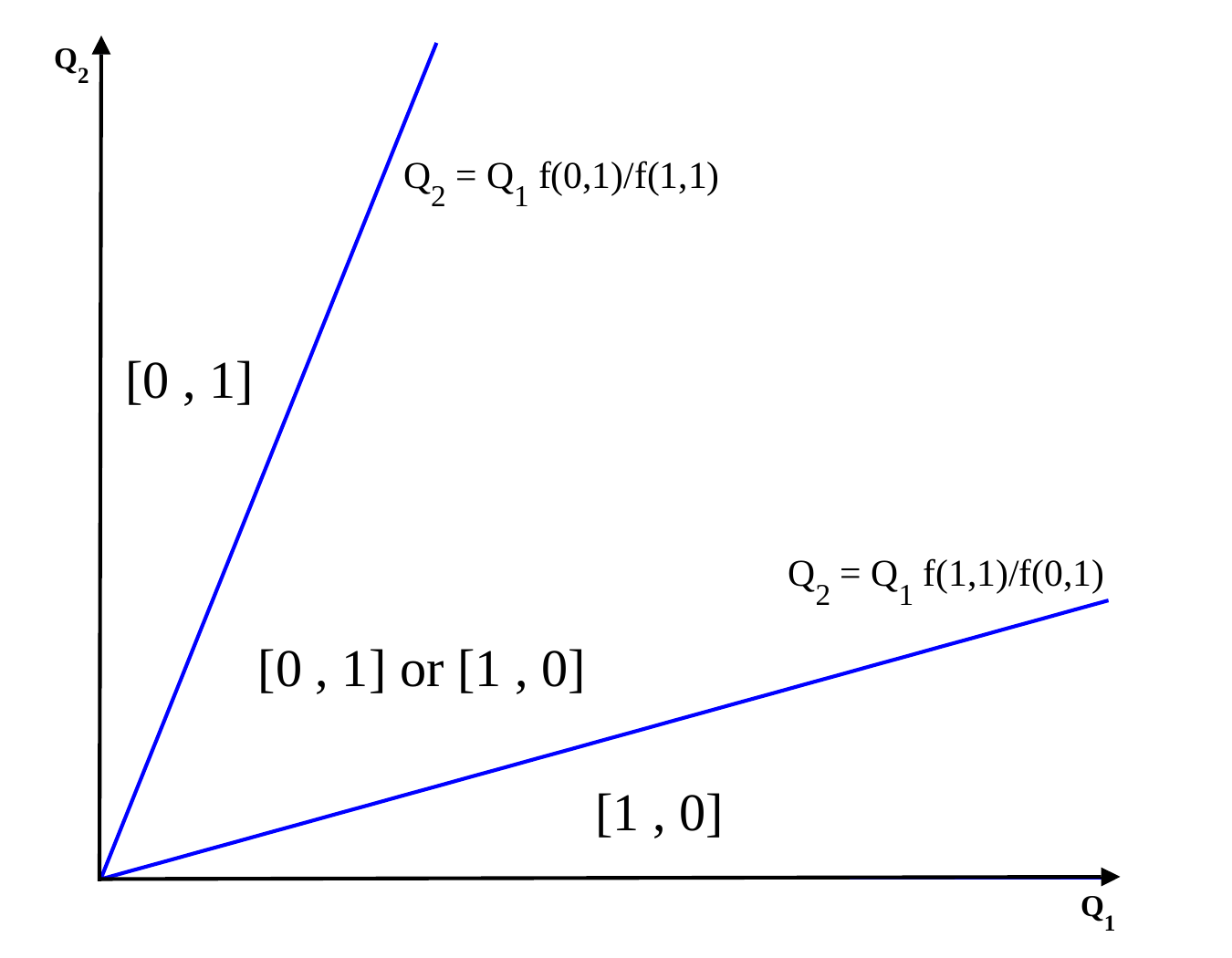}    
    \label{fig:fscheduling}
     }     
\subfigure
[LQF]
{    \includegraphics[width=0.3\textwidth]{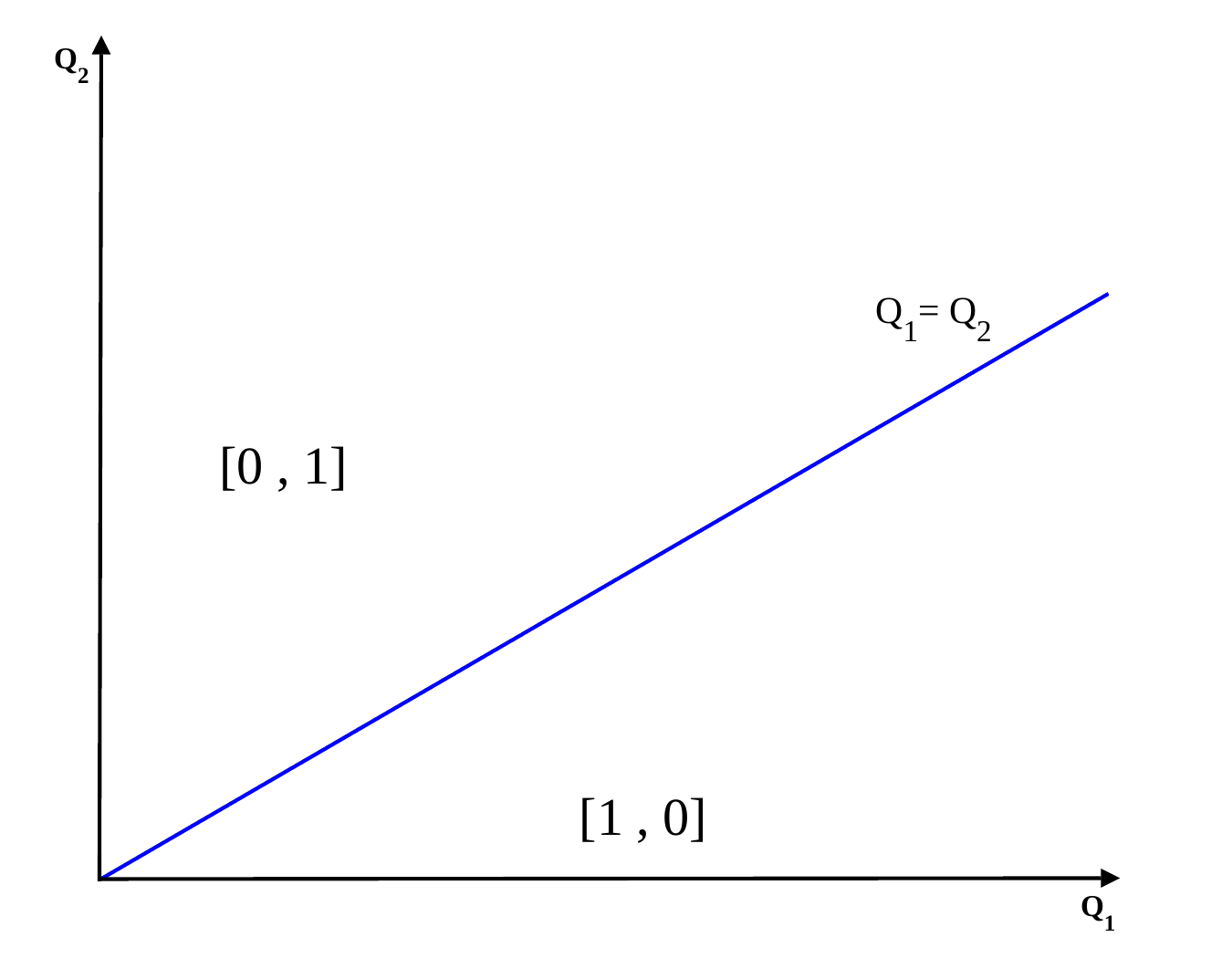}    
    \label{fig:LQF} 
}  
\subfigure
[PCS]
{    \includegraphics[width=0.3\textwidth]{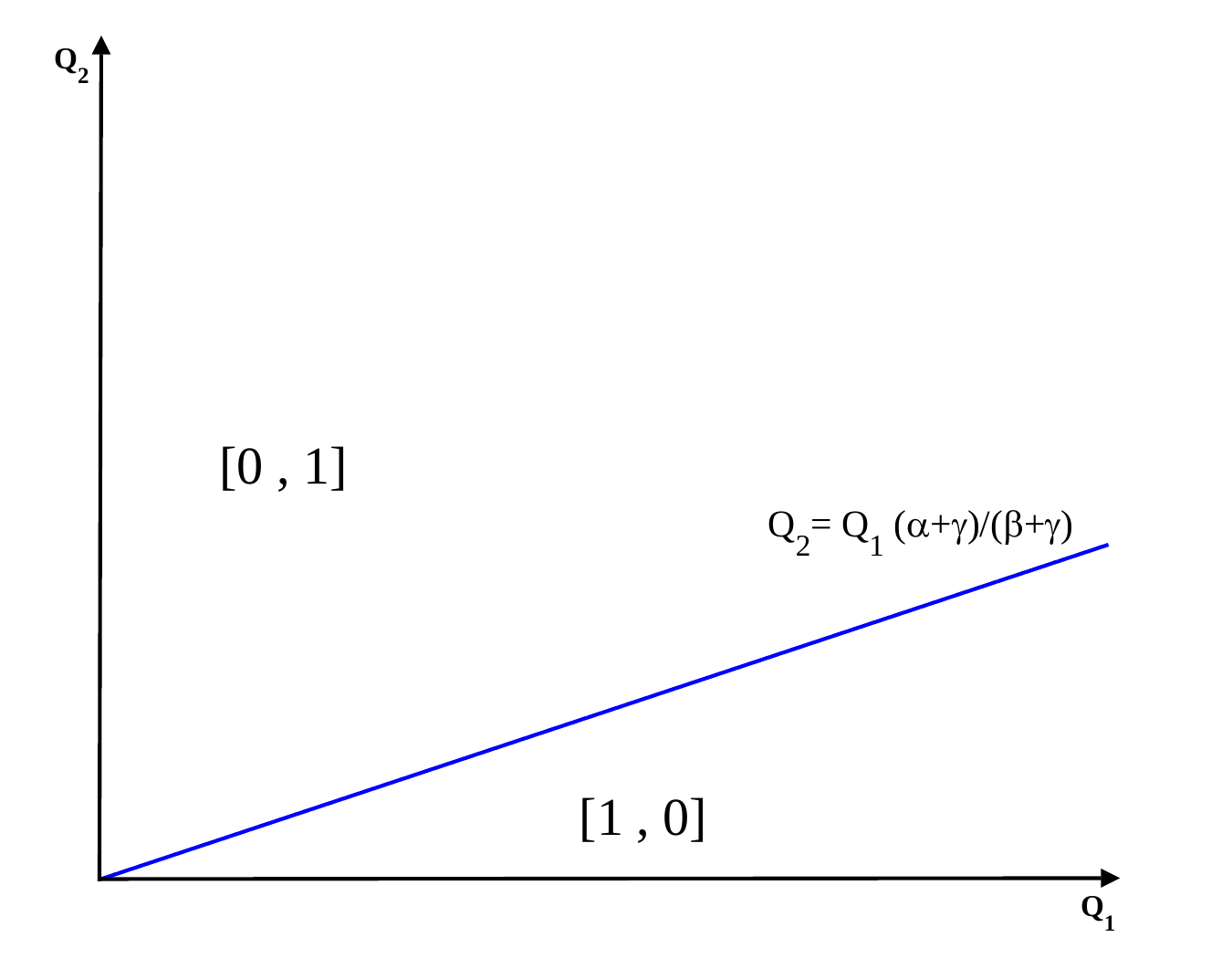}    
    \label{fig:PCS} 
}  
\caption{Scheduling decisions made by $f$-scheduling, LQF, and PCS for a $2 \times 1$ switch}
\label{fig:switch21}
\end{figure}

To compare the candidate scheduling policies, we consider a $2 \times 1$ switch.
It can be easily shown that in this case it is sufficient and necessary to consider matrix $\boldsymbol{W}$
of the form 
{\footnotesize{$\bigg [
\begin{array}{cc}
   \alpha & - \gamma  \\
   - \gamma & \beta  
\end{array} \bigg ]$}} 
where $\alpha, \beta, \gamma \ge 0$ and $\alpha \beta > \gamma^2$.  
Figure~\ref{fig:switch21} shows the scheduling decisions, $[\eta_{1}, \eta_{2}]$, made by $f$-scheduling, LQF, and PCS for different values of $\boldsymbol{Q}$ in the backlog space. 

Figure~\ref{fig:switch21} shows that for a $2 \times 1$ switch, LQF and PCS are  consistent with $f$-scheduling and for bivariate functions $f$ satisfying $\frac{f(1,1)}{f(0,1)} < \frac{\alpha+\gamma}{\beta+\gamma} < \frac{f(0,1)}{f(1,1)}$
(for instance, function $f(m,n)=\frac{1}{K^m(K^n-1)}$, $K \ge \frac{\max \{\alpha , \beta\} + \gamma}{\min \{\alpha , \beta\} + \gamma}$, is such a function). 
This result can be generalized to $M \times 1$ switches: 
for any bivariate function $f$ that satisfies conditions (C1) and (C2), LQF is always an $f$-schedule. Similarly with an 
appropriate choice of $f$, PCS is also an $f$-schedule.
In other words, Theorems~\ref{fmltdstopt} and~\ref{respectopt} provide an alternative method of proof for the throughput optimality of LQF and PCS for $M \times 1$ switches.

\begin{remarks}
For $M \times N$, $N >1$ switches, LQF and PCS policies might not, in general, be consistent with
$f$-scheduling. Moreover, the class of $f$-scheduling provides a 
new set of throughput optimal policies whose allocation 
in much of the queue state space coincide with that under a \emph{maximum size matching}\footnote{Maximum size matching is a scheduling policy under which the number of non-empty inputs that send packet to the outputs is maximized.} (note that balancing becomes necessary when the queue state visits one of the side $f$-cones with strict priorities).
\end{remarks}


\section{Discussion and Future Work}
\label{Discussion}
In this paper, we provided a large class of throughput optimal policies by considering a class of piece-wise quadratic Lyapunov functions. 
We also specialized our result to recover and prove the throughput optimality of two known routing policies, backpressure and ORCD.
The delay performance improvements of ORCD, reported in \cite{parul07}, shed light on the importance of the path-connected structure of $f$-policy.
In a parallel area of research, we have used the insight obtained by considering a path-connected $f$-policy to design throughput optimal policies with low overhead and complexity \cite{NaghshvarInfocom10,NaghshvarISCCSP10}.
For instance, an interesting research question involves the throughput and delay performance of distributed and low-complexity variants of ORCD \cite{NaghshvarInfocom10,NaghshvarISCCSP10}. 




\appendix

\subsection{Preliminary Lemmas}
\label{Preliminary}

In this appendix, we provide some preliminary lemmas. These lemmas are technical and only helpful in proving the main lemmas of the paper, i.e.\ Lemmas~\ref{existuniq}-\ref{routemax}.

\begin{lemma}
\label{lesspenalty}

Let $R=( C_1, \ldots, C_i, C_{i+1}, \ldots, C_M )$ and $R'=( C_1, \ldots, C_{i-1}, C_i \cup C_{i+1},C_{i+2}, \ldots, C_M )$ be two adjacent rank orderings. 
\begin{itemize}
	\item {If $R <_{\boldsymbol{Q}} R'$, then $$f(|C^{i-1}|,|C_i|) Q_{C_i} \le f(|C^{i-1}|,|C_i|+|C_{i+1}|) (Q_{C_i}+Q_{C_{i+1}}) \le f(|C^{i}|,|C_{i+1}|) Q_{C_{i+1}}.$$}
	\item {If $R' <_{\boldsymbol{Q}} R$, then $$f(|C^{i-1}|,|C_i|) Q_{C_i} > f(|C^{i-1}|,|C_i|+|C_{i+1}|) (Q_{C_i}+Q_{C_{i+1}}) > f(|C^{i}|,|C_{i+1}|) Q_{C_{i+1}}.$$}
\end{itemize}

\end{lemma}

\begin{IEEEproof}

Suppose $R <_{\boldsymbol{Q}} R'$. From definitions~\ref{MismatchDef} and \ref{PenalizeDef}, $m(R,R')=i$ and we have
\begin{eqnarray}
\label{lesspen1}
\Lambda_{f}(\boldsymbol{Q},R,i) \le \Lambda_{f}(\boldsymbol{Q},R',i),
\end{eqnarray} 


or equivalently
\begin{eqnarray}
\label{lesspen3}
f(|C^{i-1}|,|C_{i}|) Q_{C_i} \le f(|C^{i-1}|,|C_{i}|+|C_{i+1}|) (Q_{C_i}+Q_{C_{i+1}}) .
\end{eqnarray}

Using (\ref{lesspen3}) and property (C1) of function $f$, however,
\begin{align}
\label{lesspen4}
\nonumber
\frac{1}{f(|C^{i-1}|,|C_{i}|+|C_{i+1}|)} Q_{C_{i+1}} &=  \frac{1}{f(|C^{i-1}|,|C_{i}|} Q_{C_{i+1}} + \frac{1}{f(|C^{i-1}|+|C_{i}|,|C_{i+1}|)} Q_{C_{i+1}} \\
\nonumber
&\ge \left(\frac{1}{f(|C^{i-1}|,|C_{i}|+|C_{i+1}|)} - \frac{1}{f(|C^{i-1}|,|C_{i}|)} \right) Q_{C_{i}} + \frac{1}{f(|C^{i-1}|+|C_{i}|,|C_{i+1}|)} Q_{C_{i+1}} \\
&= \frac{1}{f(|C^{i-1}|+|C_{i}|,|C_{i+1}|)} (Q_{C_{i}}+Q_{C_{i+1}}).
\end{align}

Combining (\ref{lesspen3}) and (\ref{lesspen4}) completes the proof for the case $R <_{\boldsymbol{Q}} R'$.
Now suppose $R' <_{\boldsymbol{Q}} R$.
From Definition~\ref{PenalizeDef}, we have
\begin{eqnarray}
\label{lesspen5}
\Lambda_{f}(\boldsymbol{Q},R,i) > \Lambda_{f}(\boldsymbol{Q},R',i).
\end{eqnarray} 
The rest of the proof follows (\ref{lesspen3}) and (\ref{lesspen4}) identically.
\end{IEEEproof}


\begin{lemma}
\label{Lemma:basic1}
Let $R=(C_1,C_2, \ldots,C_M) \in \mathcal{R}_{}$ $\left (\text{or} \in \mathcal{R}^{c} \right )$ and $\boldsymbol{Q} \in D_f(R)$ $\left (\text{or} \ \boldsymbol{Q} \in D^c_f(R) \right )$. Then
\begin{eqnarray*}
f(|C^{i-1}|,|C_{i}|) Q_{C_i} \le f(|C^{i}|,|C_{i+1}|) Q_{C_{i+1}} \ \ i=1,2,\ldots,M-1.
\end{eqnarray*}
\end{lemma}

\begin{IEEEproof}
For all $1 \le i \le M-1$, $R'_i=(C_1,\ldots,C_{i-1},C_i \cup C_{i+1},C_{i+2} \ldots,C_M)$ is a one-step confinement of $R$.
Note that if $R \in \mathcal{R}^c \subseteq \mathcal{R}$ then $R'_i \in \mathcal{R}^c \subseteq \mathcal{R}$.
Now, since $\boldsymbol{Q} \in D_f(R)$ $\left ( \boldsymbol{Q} \in D^c_f(R) \right )$, we have $R <_{\boldsymbol{Q}} R'_i$ for all $1 \le i \le M-1$, and from Lemma~\ref{lesspenalty}, we have the assertion of the lemma.


\end{IEEEproof}


\begin{lemma}
\label{Lemma:basic2}
Let $R=(C_1,C_2, \ldots,C_M) \in \mathcal{R}_{}$ and $\boldsymbol{Q} \in D_f(R)$. For any node $k$ in ranking class $C_i$,
\begin{eqnarray*}
Q_k > \frac{f(0,|C^{i-1}|)}{f(|C^{i-1}|,1)} Q_{C^{i-1}} \ge Q_{C^{i-1}}.
\end{eqnarray*}\end{lemma}
 
\begin{IEEEproof}
Consider $R'=(C_1,\ldots,C_{i-1},\{k\},C_i - \{k\},C_{i+1},\ldots,C_M)$.
Since $\boldsymbol{Q} \in D_f(R)$, we have $R <_{\boldsymbol{Q}} R'$. Using Lemma~\ref{lesspenalty}, we have 
\begin{eqnarray}
\label{Lemma:basic3}
f(|C^{i-1}|,1) Q_k > f(|C^{i-1}|,|C_{i}|) Q_{C_i}.
\end{eqnarray}


On the other hand and since $\boldsymbol{Q} \in D_f(R)$, Lemma~\ref{Lemma:basic1} implies that
\begin{eqnarray}
\frac{f(|C^{i-1}|,|C_{i}|)}{f(|C^{j-1}|,|C_{j}|)} Q_{C_i} \ge Q_{C_j} \ \ j=1,2,\ldots,i-1.
\end{eqnarray}

Summing over $j=1,2,\ldots,i-1$ yields
\begin{eqnarray}
\label{basic201}
f(|C^{i-1}|,|C_{i}|) \left( \sum_{j=1}^{i-1} \frac{1}{f(|C^{j-1}|,|C_{j}|)} \right) Q_{C_i} \ge \sum_{j=1}^{i-1} Q_{C_j} = Q_{C^{i-1}}.
\end{eqnarray}

However, condition (C1) implies that
\begin{eqnarray}
\label{basic202}
\sum_{j=1}^{i-1} \frac{1}{f(|C^{j-1}|,|C_{j}|)} = \sum_{j=1}^{i-1} \frac{1}{f(\sum_{l=1}^{j-1} |C_{l}|,|C_{j}|)} = 
\frac{1}{f(0,\sum_{j=1}^{i-1} |C_{j}|)} = \frac{1}{f(0,|C^{i-1}|)}.
\end{eqnarray}

Combining (\ref{basic201}) and (\ref{basic202}), we obtain
\begin{eqnarray}
Q_{C_{i}} \ge \frac{f(0,|C^{i-1}|)}{f(|C^{i-1}|,|C_{i}|)} Q_{C^{i-1}} 
, 
\end{eqnarray}
which together with (\ref{Lemma:basic3}) and condition (C2) completes the proof. 

\end{IEEEproof}


\begin{lemma}
\label{Lemma:basic4}
Let $R=(C_1,C_2, \ldots,C_M) \in \mathcal{R}^{c}$ and $\boldsymbol{Q} \in D^c_f(R)$. For any node $k$ in ranking class $C_i$ that reaches 
a node in $C^{i-1}$,
\begin{eqnarray*}
Q_k > \frac{f(0,|C^{i-1}|)}{f(|C^{i-1}|,1)} Q_{C^{i-1}} \ge Q_{C^{i-1}}.
\end{eqnarray*}\end{lemma}
 
\begin{IEEEproof}
Consider $R'=(C_1,\ldots,C_{i-1},\{k\},C_i - \{k\},C_{i+1},\ldots,C_M)$. Note that $R'$ is path-connected since $R$ is path-connected and node $k$ reaches a node in $C^{i-1}$.  
Since $\boldsymbol{Q} \in D^c_f(R)$, we have $R <_{\boldsymbol{Q}} R'$, which together with Lemma~\ref{lesspenalty} gives, 
\begin{eqnarray}
\label{Lemma:basic5}
f(|C^{i-1}|,1) Q_k > f(|C^{i-1}|,|C_{i}|) Q_{C_i}.
\end{eqnarray}


The rest of the proof is similar to the proof of Lemma~\ref{Lemma:basic2} and is omitted for brevity.

\end{IEEEproof}


\begin{lemma}
\label{existuniq_Claim2}
Let $R=(C_1,\ldots,C_i,C_{i+1},\ldots,C_M)$ and $\tilde{R}=(C_1,\ldots,C_{i-1},C_i \cup C_{i+1},C_{i+2},\ldots,C_M)$. Suppose following assumptions hold:
\begin{enumerate}
	\item { $R <_{\boldsymbol{Q}} R'$ for all $R' \in \mathcal{B}_1(R)$.}
	\item {
	For any node $k \in C_{i+1}$, $Q_k > Q_{C_i}$.}
	\item {$\tilde{R} <_{\boldsymbol{Q}} R$.}
\end{enumerate}
Then $\tilde{R} <_{\boldsymbol{Q}} R'$ for all $R' \in \mathcal{B}_1(\tilde{R})$.
\end{lemma}

\begin{IEEEproof}

It is sufficient to show that $\tilde{R}$ penalizes $\boldsymbol{Q}$ less than its one-step refinements with regard to $C_i \cup C_{i+1}$.
Let $\tilde{R}_1=(C_1,\ldots,C_{i-1},A \cup C,B \cup D,C_{i+2},\ldots,C_M)$ be a one-step refinement of $\tilde{R}$ where $A$, $B$, $C$, $D$ are sets of nodes satisfying $A \neq \emptyset$, $D \neq \emptyset$, $C_i=A \cup B$, and $C_{i+1}=C \cup D$. Then we can write $R$ and $\tilde{R}$ as $R=(C_1,\ldots,C_{i-1},A \cup B,C \cup D,C_{i+2},\ldots,C_M)$ and $\tilde{R}=(C_1,\ldots,C_{i-1},A \cup B \cup C \cup D,C_{i+2},\ldots,C_M)$. Let $R_1=(C_1,\ldots,C_{i-1},A,B,C \cup D,C_{i+2},\ldots,C_M)$ and $R_2=(C_1,\ldots,C_{i-1},A \cup B,C,D,C_{i+2},\ldots,C_M)$ be one-step refinements of $R$. 
Let $\sum_{j=1}^{i-1}|C_{j}|=m$, $|A|=a$, $|B|=b$, $|C|=c$, and $|D|=d$. 
We consider three cases based on sets $B$ and $C$:

{\bf{Case I.}} $B$ and $C$ are not empty.

Since $R <_{\boldsymbol{Q}} R_2$, by lemma \ref{lesspenalty}, we have 
\begin{eqnarray}
\label{existuniq_Lemma2_01}
f(m+a+b,c) Q_C > f(m+a+b,c+d) (Q_C + Q_D) > f(m+a+b+c,d)  Q_D.
\end{eqnarray}

Let assume that $\tilde{R}_1 <_{\boldsymbol{Q}} \tilde{R}$. By Lemma~\ref{lesspenalty}, we have
\begin{eqnarray}
\label{existuniq_Lemma2_02}
f(m,a+c) (Q_A + Q_C) \le f(m+a+c,b+d) (Q_B + Q_D).
\end{eqnarray}

After proper arrangement,
\begin{eqnarray}
\label{existuniq_Lemma2_03}
Q_B \ge \frac{f(m,a+c)}{f(m+a+c,b+d)} Q_C - Q_D.
\end{eqnarray}

By property (C1) of function $f$,
\begin{eqnarray}
\label{existuniq_Lemma2_04}
\frac{1}{f(m+a+c,b+d)} = \frac{1}{f(m+a+c,b)} + \frac{1}{f(m+a+b+c,d)}.
\end{eqnarray}

Combining (\ref{existuniq_Lemma2_01}), (\ref{existuniq_Lemma2_03}), and (\ref{existuniq_Lemma2_04}), we obtain 
\begin{eqnarray}
\label{existuniq_Lemma2_05}
\nonumber
Q_B &>& \left( \frac{f(m,a+c)}{f(m+a+c,b+d)} - \frac{f(m+a+b,c)}{f(m+a+b+c,d)} \right) Q_C \\
&=& \left( \frac{f(m,a+c)}{f(m+a+c,b)} + \frac{f(m,a+c)-f(m+a+b,c)}{f(m+a+b+c,d)} \right) Q_C.
\end{eqnarray}

By property (C2) of function $f$, 
\begin{eqnarray}
\label{existuniq_Lemma2_06}
\frac{f(m,a+c)}{f(m+a+c,b)} \ge 1.
\end{eqnarray}

By Property (C1) and (C2) of function $f$,
\begin{eqnarray}
\label{existuniq_Lemma2_07}
\nonumber
\frac{2}{f(m,a+c)} &\le& \frac{1}{f(m,a+c)} + \frac{1}{f(m+a+c,b)} \\
\nonumber
&=& \frac{1}{f(m,a+c+b)} \\ 
&=& \frac{1}{f(m,a+b)} + \frac{1}{f(m+a+b,c)} \le \frac{2}{f(m+a+b,c)}.
\end{eqnarray}

Combining (\ref{existuniq_Lemma2_05}), (\ref{existuniq_Lemma2_06}), and (\ref{existuniq_Lemma2_07}), we obtain 
\begin{eqnarray}
\label{existuniq_Lemma2_08}
Q_B > Q_C. 
\end{eqnarray}

By assumption 2 of the lemma, queue backlog of any node in set $C$ is larger than $Q_A+Q_B$. But this is in contradiction with (\ref{existuniq_Lemma2_08}). Therefore, assumption $\tilde{R}_1 <_{\boldsymbol{Q}} \tilde{R}$ cannot hold and we have $\tilde{R} <_{\boldsymbol{Q}} \tilde{R}_1$.  

{\bf{Case II.}} $B$ is empty.

Since $\tilde{R} <_{\boldsymbol{Q}} R$, we have the following inequality by Lemma~\ref{lesspenalty},
\begin{eqnarray}
\label{existuniq_Lemma2_09}
f(m,a+b) (Q_A + Q_B) > f(m+a+b,c+d) (Q_C + Q_D).
\end{eqnarray}

Using (\ref{existuniq_Lemma2_01}), (\ref{existuniq_Lemma2_09}), and the fact that $B=\emptyset$, we obtain following inequalities
\begin{eqnarray}
\label{existuniq_Lemma2_10}
f(m+a,c)Q_C &>& f(m+a+c,d) Q_D, \\
\label{existuniq_Lemma2_11}
f(m,a) Q_A &>& f(m+a+c,d) Q_D.
\end{eqnarray}

By (\ref{existuniq_Lemma2_10}), (\ref{existuniq_Lemma2_11}), and property $(C1)$ of function $f$,
\begin{eqnarray}
\label{existuniq_Lemma2_12}
\nonumber
f(m,a+c) (Q_A+Q_C) &>& f(m,a+c) \left( \frac{1}{f(m,a)} + \frac{1}{f(m+a,c)} \right) f(m+a+c,d) Q_D \\
&=& f(m+a+c,d) Q_D.
\end{eqnarray}

By Lemma \ref{lesspenalty}, $\tilde{R} <_{\boldsymbol{Q}} \tilde{R}_1$.

{\bf{Case III.}} $C$ is empty.

Since $R <_{\boldsymbol{Q}} R_1$, we have
\begin{eqnarray}
\label{existuniq_Lemma2_13}
f(m,a) Q_A > f(m,a+b) (Q_A + Q_B) > f(m+a,b) Q_B.
\end{eqnarray}

Using (\ref{existuniq_Lemma2_02}), (\ref{existuniq_Lemma2_13}), and the fact that $C=\emptyset$ we obtain following inequalities
\begin{eqnarray}
\label{existuniq_Lemma2_14}
f(m,a) Q_A &>& f(m+a,b) Q_B, \\
\label{existuniq_Lemma2_15}
f(m,a) Q_A &>& f(m+a+b,d) Q_D. 
\end{eqnarray}

Combining (\ref{existuniq_Lemma2_14}) and (\ref{existuniq_Lemma2_15}), we obtain
\begin{eqnarray}
\label{existuniq_Lemma2_16}
\nonumber
Q_B + Q_D &<& f(m,a) \left( \frac{1}{f(m+a,b)} + \frac{1}{f(m+a+b,d)} \right) Q_A \\
&=& \frac{f(m,a)}{f(m+a,b+d)} Q_A.
\end{eqnarray}

By Lemma~\ref{lesspenalty}, $\tilde{R} <_{\boldsymbol{Q}} \tilde{R}_1$.

\end{IEEEproof}


\begin{lemma}
\label{existuniq_Claim2new}
Let $R=(C_1,\ldots,C_i,C_{i+1},\ldots,C_M)\in \mathcal{R}^c$ and $\tilde{R}=(C_1,\ldots,C_{i-1},C_i \cup C_{i+1},C_{i+2},\ldots,C_M)$. Suppose following assumptions hold:
\begin{enumerate}
	\item { $R <_{\boldsymbol{Q}} R'$ for all $R' \in \mathcal{B}_1^c(R)$.}
	\item {For any node $k \in C_{i+1}$ such that it reaches a node in $C^i$, $Q_k > Q_{C_i}$.}
	\item {$\tilde{R} <_{\boldsymbol{Q}} R$.}
\end{enumerate}
Then $\tilde{R} <_{\boldsymbol{Q}} R'$ for all $R' \in \mathcal{B}_1^c(\tilde{R})$.
\end{lemma}
 
The proof of Lemma~\ref{existuniq_Claim2new} is very similar to the proof of Lemma~\ref{existuniq_Claim2}.
Here we let $\tilde{R}_1=(C_1,\ldots,C_{i-1},A \cup C,B \cup D,C_{i+2},\ldots,C_M)$ be a path-connected one-step refinement of $\tilde{R}$ where $A$, $B$, $C$, $D$ are sets of nodes satisfying $A \neq \emptyset$, $D \neq \emptyset$, $C_i=A \cup B$, and $C_{i+1}=C \cup D$. 
Let $R_1=(C_1,\ldots,C_{i-1},A,B,C \cup D,C_{i+2},\ldots,C_M)$ and $R_2=(C_1,\ldots,C_{i-1},A \cup B,C,D,C_{i+2},\ldots,C_M)$ be one-step refinements of $R$. We know that $R$ and $\tilde{R}_1$ are path-connected. This implies that $R_2$ is path-connected and $R_1$ is path-connected when $C=\emptyset$.
To get a contradiction to (\ref{existuniq_Lemma2_08}), we should note that since $R_2$ is path-connected, there exists at least one node in $C$ such that it reaches a node in $C^i$ and hence, by assumption 2 of Lemma~\ref{existuniq_Claim2new} we have $Q_C > Q_A + Q_B$.


\subsection{Proof of Lemmas~\ref{existuniq} and \ref{existuniq_critical}}
\label{App:existuniq}

This appendix is dedicated to the proof of Lemmas~\ref{existuniq} and \ref{existuniq_critical}. These lemmas contain extended algebraic manipulation to show that cones $D_f(\cdot)$ partition $\mathbb{R}^N_{+}$, i.e. for $\forall \boldsymbol{Q} \in \mathbb{R}^N_{+}$, $\exists! \ R \in \mathcal{R}$ such that $\boldsymbol{Q} \in D_f(R)$. 
The existence proofs are inductive, while the uniqueness proofs are done by contradiction. 

\hspace{-0.2 in} {\bf{Lemma~\ref{existuniq}}}. {\textit{
Let bivariate function $f$ satisfy conditions (C1) and (C2).
Then for all $\boldsymbol{Q} \in \mathbb{R}^{N}_{+}$, there exists a unique $R \in \mathcal{R}_{}$ such that $\boldsymbol{Q} \in D_f(R)$.
}}

{\bf{Proof of Existence:}}

The proof is done by induction.
Let $n$ denote the total number of nodes in the network excluding the destination.
For $n=1$ there exists only one rank ordering and the proof for this case is trivial. Now suppose for all $n \le N-1$ and all $\boldsymbol{Q} \in \mathbb{R}^{n}_+$ there exists a rank ordering $R$ such that $\boldsymbol{Q} \in D_f(R)$. 
Next we constructively show that for all $\boldsymbol{Q} \in \mathbb{R}^{N}_{+}$, there exists a rank ordering $R$ such that $\boldsymbol{Q} \in D_f(R)$.

\begin{enumerate}
	\item[1.] {Let $R_0=( \{ 1,2,\ldots,N \} )$.}
	\item[2.] {
\begin{enumerate}
  \item[]{}
	\item[2.1.] {Initialize $l=1$.}
	\item[2.2.] {Is there a rank ordering $\hat{R}$ of the form $\hat{R}=(\hat{C}_1,\hat{C}_2)$, where $|\hat{C}_1|=N-l$, $|\hat{C}_2|=l$, and $\hat{R} <_{\boldsymbol{Q}} R_0$?}
	\item[2.3.] {If yes, go to step 3. Otherwise, go to step 2.4.}
	\item[2.4.] {$l=l+1$. Is $l < N$?}
	\item[2.5.] {If yes, go to step 2.2. Otherwise, $\boldsymbol{Q} \in D_f(R_0)$.}
\end{enumerate}
	\item[3.] {Consider nodes in class $\hat{C}_1$ of rank ordering $\hat{R}$. Since $|\hat{C}_1|<N$, by the assumption of the induction, there exists a rank ordering for the nodes in $\hat{C}_1$ such that it penalizes $\boldsymbol{Q}$ less than all its adjacent rank orderings. Let $R^*=(C^*_1,C^*_2,\ldots,C^*_{M-1})$ be this rank ordering. Let $R^*_0=(C^*_1,C^*_2,\ldots,C^*_M)=(C^*_1,C^*_2,\ldots,C^*_{M-1},\hat{C}_2)$. 
Furthermore, let $R^*_i=(C_1,C_2,\ldots,C_{M-i-1},C_{M-i} \cup \ldots \cup C_M)$ denote the rank ordering generated by merging the last $i$ classes of $R^*_0$.}
	\item[4.] { Find $m$ such that $R^*_i <_{\boldsymbol{Q}} R^*_{i-1}$ for for $i=1,2, \ldots, m$, but $R^*_m <_{\boldsymbol{Q}} R^*_{m+1}$.
	Claim~\ref{existuniq_Claim1} below and Lemma~\ref{existuniq_Claim2} in Appendix~\ref{Preliminary} establish that $\boldsymbol{Q} \in D_f(R^*_m)$.}
}
\end{enumerate}
   
\begin{claim}
\label{existuniq_Claim1}
$R^*_0 <_{\boldsymbol{Q}} R'$ for all $R' \in \mathcal{B}_1(R^*_0)$.
Moreover, for $i=1,2,\ldots,M-1$, and for any node $k \in C^*_{i+1}$, $Q_k > \sum_{j=1}^{i} Q_{C^*_{j}}$.
\end{claim}   

By Claim~\ref{existuniq_Claim1}, Lemma~\ref{existuniq_Claim2}, and using the fact that $R^*_i <_{\boldsymbol{Q}} R^*_{i-1}$ \ $i=1,2, \ldots, m$, we can recursively show that $R^*_i <_{\boldsymbol{Q}} R'$ for all $R' \in \mathcal{B}_1(R^*_i)$, for $i=1,2, \ldots, m$.
By construction, we also know that $R^*_m <_{\boldsymbol{Q}} R^*_{m+1}$.
Moreover, $R^*_m$ penalizes $\boldsymbol{Q}$ less than its one-step confinements with regard to $C^*_i$, $i=1,2,\ldots,m-2$, since $R^* <_{\boldsymbol{Q}} R'$ for all $R' \in \mathcal{B}_2(R^*)$.
Hence, $R^*_m <_{\boldsymbol{Q}} R'$ for all $R' \in \mathcal{A}(R^*_m)$, and by definition, $\boldsymbol{Q} \in D_f(R^*_m)$.
Now what remains is to verify Claim~\ref{existuniq_Claim1}.



\begin{IEEEproof}[Proof of Claim \ref{existuniq_Claim1}]
Note that following results are immediate using Lemmas \ref{Lemma:basic1}, \ref{Lemma:basic2}, and the fact that 
$R^* <_{\boldsymbol{Q}} R'$ for all $R' \in \mathcal{B}_1(R^*)$:
\begin{itemize}
	\item {$R^*_0$ penalizes $\boldsymbol{Q}$ less than all its one-step refinements with regard to ranking class $C^*_i$ for $i=1,2,\ldots,M-1$.}
	\item {For $i=1,2,\ldots,M-2$, and for any node $k \in C^*_{i+1}$, $Q_k > \sum_{j=1}^{i} Q_{C^*_{j}}$ .}
\end{itemize}

What is left is to show that
\begin{enumerate}
	\item {$R^*_0$ penalizes $\boldsymbol{Q}$ less than all its one-step refinements with regard to ranking class $C^*_M$.}
	\item {For any node $k \in C^*_{M}$, $Q_k > \sum_{j=1}^{M-1} Q_{C^*_{j}}$.}
\end{enumerate}

Let $\tilde{R}=(C^*_1,C^*_2,\ldots,C^*_{M-1},A,B)$ be a one-step refinement of $R^*_0$ with regard to $C^*_M$, i.e. $A \cup B = C^*_M$.
Note that $\cup_{i=1}^{M-1} C^*_i = \hat{C}_1$ and $C^*_M = \hat{C}_2$. 
Suppose $\tilde{R} <_{\boldsymbol{Q}} R^*_0$. By Lemma~\ref{lesspenalty}, we have
\begin{align}
\label{existuniq_Lemma1_01}
	f(|\hat{C}_1|,|A|) Q_{A} \le f(|\hat{C}_1|,|\hat{C}_2|) Q_{\hat{C}_2} \le f(|\hat{C}_1| + |A|,|B|) Q_{B}.
\end{align}

On the other hand, since $\hat{R} <_{\boldsymbol{Q}} R_0$, by Lemma~\ref{lesspenalty}, we have
\begin{align}
\label{existuniq_Lemma1_02}
	f(0,|\hat{C}_1|) Q_{\hat{C}_1} \le f(|\hat{C}_1|,|\hat{C}_2|) Q_{\hat{C}_2}.
\end{align}

Combining (\ref{existuniq_Lemma1_01}) and (\ref{existuniq_Lemma1_02}), and using property (C1) of function $f$, we obtain
\begin{align}
\label{existuniq_Lemma1_03}
\nonumber
\frac{f(|\hat{C}_1| + |A|,|B|)}{f(0,|\hat{C}_1|+|A|)} Q_{B} &=
f(|\hat{C}_1| + |A|,|B|) Q_{B} \left ( \frac{1}{f(0,|\hat{C}_1|)} + \frac{1}{f(|\hat{C}_1|,|A|)} \right )	\\
& \ge Q_{\hat{C}_1} + Q_{A}.
\end{align}

After proper arrangement we have 
\begin{align}
\label{existuniq_Lemma1_04}
f(0,|\hat{C}_1|+|A|) \left ( Q_{\hat{C}_1} + Q_{A} \right )  \le f(|\hat{C}_1| + |A|,|B|) Q_{B},
\end{align}
which implies that rank ordering $(\hat{C}_1 \cup A, B)$ penalizes $\boldsymbol{Q}$ less than $R_0$. But this is a contradiction (look at step 2 of the given procedure and note that $|B| < |\hat{C}_2|$ ). 
Therefore, $R^*_0$ penalizes $\boldsymbol{Q}$ less than all its one-step refinements with regard to ranking class $C^*_M$.

Now consider node $k \in C^*_M$ and let $\tilde{R}=(C^*_1,C^*_2,\ldots,C^*_{M-1},\{ k \},C^*_{M} - \{ k \})$. From the result of the previous part, $R^*_0 <_{\boldsymbol{Q}} \tilde{R}$. By Lemma~\ref{lesspenalty}, we have
\begin{align}
\label{existuniq_Lemma1_05}
	f(|\hat{C}_1|,1) Q_{k} > f(|\hat{C}_1|,|\hat{C}_2|) Q_{\hat{C}_2}.
\end{align}
 
Combining (\ref{existuniq_Lemma1_02}) and (\ref{existuniq_Lemma1_05}), we have
\begin{align}
\label{existuniq_Lemma1_06}
	 Q_{k} > \frac{f(0,|\hat{C}_1|)}{f(|\hat{C}_1|,1)} Q_{\hat{C}_1} \ge Q_{\hat{C}_1},
\end{align}
where the last inequality follows from property (C2) of function $f$.
Hence, for all  $k \in C^*_{M}$, $Q_k > \sum_{j=1}^{M-1} Q_{C^*_{j}}$.

\end{IEEEproof}

{\bf{Proof of Uniqueness:}}

Consider $R=\{ C_1, C_2, \ldots, C_M \}$ and $\hat{R}=\{ \hat{C}_1, \hat{C}_2, \ldots, \hat{C}_{\hat{M}} \}$. We will prove by contradiction that $\boldsymbol{Q}$ cannot be in $D_f(R)$ and $D_f(\hat{R})$ simultaneously.

{\bf{Case I.}} There exist nodes $a$ and $b$ such that $b \prec^{R} a$ and $a \prec^{\hat{R}} b$, i.e.
\[ \left\{ \begin{array}{ll}
      b \in C^{i-1} &, \ a \in \cup_{l=i}^{M} C_{l} \\
      a \in \hat{C}^{j-1} &, \ b \in \cup_{l=j}^{\hat{M}} \hat{C}_{l}          
         \end{array} \right. . \]

If $\boldsymbol{Q} \in D_f(R)$, by Lemma~\ref{Lemma:basic2}, we have
\begin{eqnarray}
\label{uniquenessA1}
Q_a > Q_{C^{i-1}} \ge Q_b.
\end{eqnarray}

Similarly, if $\boldsymbol{Q} \in D_f(\hat{R})$, by Lemma~\ref{Lemma:basic2}, we have
\begin{eqnarray}
\label{uniquenessA2}
Q_b > Q_{\hat{C}^{j-1}} \ge Q_a.
\end{eqnarray}

Clearly (\ref{uniquenessA1}) and (\ref{uniquenessA2}) cannot hold simultaneously. 

{\bf{Case II.}} There are no nodes $a$, $b$, such that $b \prec^{R} a$ and $a \prec^{\hat{R}} b$. In this case, it is not difficult to see that, there exist $n$, $n \le M$, consecutive classes $C_{i+1}, \ldots , C_{i+n} \in R$, and $\hat{C}_j \in \hat{R}$ such that for some sets of nodes $A_1$, $A_2$, $B_1$, $\ldots$, $B_n$, the following relationships hold
\[ \left\{ \begin{array}{ll}
      C_{i+1} &= A_1 \cup B_1 \\
      C_{i+l} &= B_l \ \ \ 2 \le l \le n-1 \\ 
      C_{i+n} &= B_n \cup A_2         
         \end{array} \right. , \]
and
\begin{eqnarray}
\nonumber
\hat{C}_j= \cup_{l=1}^{n} B_l,
\end{eqnarray} 
where $B_1$, $\ldots$, $B_n$ are non-empty while $A_1$ and $A_2$ could be empty.

In rank ordering $R$, $C^{i}$ and $A_1$ have lower rank than $\cup_{l=2}^{n} B_l$. 
Because of the condition of Case II, none of the nodes in $C^{i} \cup A_1$ can have a higher rank than a node in $\cup_{l=2}^{n} B_l$ under rank ordering $\hat{R}$. Hence, we have 
\begin{eqnarray}
C^{i} \cup A_1 &=& \hat{C}^{j-1}, \\
\label{uniquenessB11}
|C^{i}| + |A_1| &=& |\hat{C}^{j-1}|.
\end{eqnarray}

Furthermore,
\begin{eqnarray}
\label{uniquenessB12}
\nonumber
|C^{i+n-1}| &=& |C^{i}| + |A_1| + |\cup_{l=1}^{n-1}B_l| \\
&=&  |\hat{C}^{j-1}| + \sum_{l=1}^{n-1} |B_l|.
\end{eqnarray}

Now suppose $\boldsymbol{Q} \in D_f(R)$.
Let $R_1= \{ C_1, \ldots, C_{i}, A_1, B_1, C_{i+2}, \ldots, C_M \}$ and $R_2= \{ C_1 , \ldots , C_{i+n-1} , B_n , A_2 , C_{i+n+1} ,$ $\ldots , C_M \}$ be one-step refinements of $R$. 
Since $\boldsymbol{Q} \in D_f(R)$, $R <_{\boldsymbol{Q}} R_1$ and $R <_{\boldsymbol{Q}} R_2$. By Lemma~\ref{lesspenalty}, we have
\begin{eqnarray}
\label{uniquenessB1}
f(|C^{i}|,|C_{i+1}|) Q_{C_{i+1}} &\ge& f(|C^{i}|+|A_1|,|B_1|) Q_{B_1}, \\
\label{uniquenessB2}
f(|C^{i+n-1}|,|B_n|) Q_{B_n} &\ge& f(|C^{i+n-1}|,|C_{i+n}|) Q_{C_{i+n}}, 
\end{eqnarray}
where equality in (\ref{uniquenessB1}) and (\ref{uniquenessB2}) hold when $A_1$ and $A_2$ are empty respectively. 
Moreover, since $\boldsymbol{Q} \in D_f(R)$, by Lemma~\ref{Lemma:basic1} we have 
\begin{eqnarray}
\label{uniquenessB3}
f(|C^{l}|,|C_{l+1}|) Q_{C_{l+1}} \ge f(|C^{l-1}|,|C_l|) Q_{C_l}  \ \ l=1,2,\ldots,M-1. 
\end{eqnarray}

Combining (\ref{uniquenessB1})-(\ref{uniquenessB3}), we obtain
\begin{eqnarray}
\label{uniquenessB4}
f(|C^{i+n-1}|,|B_n|) Q_{B_n} \ge f(|C^{i}|+|A_1|,|B_1|) Q_{B_1}.
\end{eqnarray} 

However, we also have assumed that $\boldsymbol{Q} \in D_f(\hat{R})$. 
Let $\hat{R}_1=\{ \hat{C}_1, \ldots, \hat{C}_{j-1}, B_1, \cup_{l=2}^{n} B_l, \hat{C}_{j+1}, \ldots, \hat{C}_{\hat{M}} \}$ and 
{$\hat{R}_2=\{ \hat{C}_1, \ldots, \hat{C}_{j-1}, \cup_{l=1}^{n-1} B_l, B_n, \hat{C}_{j+1}, \ldots, \hat{C}_{M'} \}$} be one-step refinements of $\hat{R}$. By Lemma~\ref{lesspenalty}, we have 
\begin{eqnarray}
\label{uniquenessB5}
f(|\hat{C}^{j-1}|,\sum_{l=1}^{n}|B_l|) \sum_{l=1}^{n} Q_{B_l} &<& f(|\hat{C}^{j-1}|,|B_1|) Q_{B_1}, \\
\label{uniquenessB6}
f(|\hat{C}^{j-1}|+\sum_{l=1}^{n-1}|B_l|,|B_n|) Q_{B_n} &<& f(|\hat{C}^{j-1}|,\sum_{l=1}^{n}|B_l|) \sum_{l=1}^{n} Q_{B_l},
\end{eqnarray}
whose direct consequence is
\begin{eqnarray}
\label{uniquenessB7}
f(|\hat{C}^{j-1}|+\sum_{l=1}^{n-1}|B_l|,|B_n|) Q_{B_n} < f(|\hat{C}^{j-1}|,|B_1|) Q_{B_1}.
\end{eqnarray}
Substituting (\ref{uniquenessB11}) and (\ref{uniquenessB12}) in  (\ref{uniquenessB7}), we obtain
\begin{eqnarray}
\label{uniquenessB9}
f(|C^{i+n-1}|,|B_n|) Q_{B_n} < f(|C^{i-1}|+|A_1|,|B_1|) Q_{B_1},
\end{eqnarray}
which contradicts (\ref{uniquenessB4}). 
Therefore, $\boldsymbol{Q}$ cannot be in $D_f(R)$ and $D_f(\hat{R})$ simultaneously.



Next, we provide a brief sketch of the proof of Lemma~\ref{existuniq_critical}.

\hspace{-0.2 in} {\bf{Lemma~\ref{existuniq_critical}}}. {\textit{
If bivariate function $f$ satisfies conditions (C1) and (C2), then for all $\boldsymbol{Q} \in \mathbb{R}^{N}_{+}$, there exists a unique $R \in \mathcal{R}^{c}$ such that $\boldsymbol{Q} \in D^c_f(R)$.
}}

Note that this proof is very similar to the proof of Lemma~\ref{existuniq} and the only difference is the limitation to the set of path-connected rank orderings. To prove Lemma~\ref{existuniq_critical}, we only need to show that the rank orderings used in Lemma~\ref{existuniq} can be identically selected from path-connected rank orderings.


{\bf{Proof of Existence:}}

The proof is done by induction.
Let $n$ denote the total number of nodes in the network excluding the destination.
For $n=1$ there exists only one rank ordering and the proof for this case is trivial. Now suppose for $n \le N-1$ and for any $\boldsymbol{Q} \in \mathbb{R}^{n}_+$ there exists a path-connected rank ordering $R$ such that $R <_{\boldsymbol{Q}} R'$ for all $R' \in \mathcal{A}^{c}(R)$. 
For $n=N$ and for any $\boldsymbol{Q} \in \mathbb{R}^{N}_{+}$, using the procedure below, we will constructively show that there exists a path-connected rank ordering $R$ such that $\boldsymbol{Q} \in D_f^c(R)$. 

\begin{enumerate}
	\item[1.] {Let $R_0=( \{ 1,2,\ldots,N \} )$.}
	\item[2.] {
\begin{enumerate}
  \item[]{}
	\item[2.1.] {Initialize $l=1$.}
	\item[2.2.] {Is there a path-connected rank ordering $\hat{R}$ of the form $\hat{R}=(\hat{C}_1,\hat{C}_2)$, where $|\hat{C}_1|=N-l$, $|\hat{C}_2|=l$, and $\hat{R} <_{\boldsymbol{Q}} R_0$?}
	\item[2.3.] {If yes, go to step 3. Otherwise, go to step 2.4.}
	\item[2.4.] {$l=l+1$. Is $l < N$?}
	\item[2.5.] {If yes, go to step 2.2. Otherwise, $\boldsymbol{Q} \in D_f^c (R_0)$.}
\end{enumerate} }
	\item[3.] {Consider nodes in class $\hat{C}_1$ of rank ordering $\hat{R}$. Since $|\hat{C}_1|<N$, by the assumption of the induction, there exists a path-connected rank ordering for the nodes in $\hat{C}_1$ such that it penalizes $\boldsymbol{Q}$ less than all its adjacent path-connected rank orderings. Let $R^*=(C^*_1,C^*_2,\ldots,C^*_{M-1})$ be this rank ordering. Let $R^*_0=(C^*_1,C^*_2,\ldots,C^*_M)=(C^*_1,C^*_2,\ldots,C^*_{M-1},\hat{C}_2)$. 
Furthermore, let $R^*_i=(C_1,C_2,\ldots,C_{M-i-1},C_{M-i} \cup \ldots \cup C_M)$ denote the rank ordering generated by merging the last $i$ classes of $R^*_0$. (It is clear that rank orderings $R^*_i$, $i=0,1,\ldots,M-1$, are path-connected.)}
	\item[4.] { Find $m$ such that $R^*_i <_{\boldsymbol{Q}} R^*_{i-1}$ for for $i=1,2, \ldots, m$, but $R^*_m <_{\boldsymbol{Q}} R^*_{m+1}$.
	Modified version of Claim~\ref{existuniq_Claim1} below and Lemma~\ref{existuniq_Claim2new} in Appendix~\ref{Preliminary} establish that $\boldsymbol{Q} \in D_f^c (R^*_m)$. }
\end{enumerate}
   
\begin{claim}
\label{existuniq_Claim1new} 
$R^*_0 <_{\boldsymbol{Q}} R'$ for all $R' \in \mathcal{B}_1^c(R^*_0)$.
Moreover, for $i=1,2,\ldots,M-1$, and for any node $k \in C^*_{i+1}$ such that it reaches a node in $\cup_{j=1}^{i} C_j^*$, $Q_k > \sum_{j=1}^{i} Q_{C^*_{j}}$.
\end{claim}

The proof of Claim~\ref{existuniq_Claim1new} is very similar to the proof of Claim~\ref{existuniq_Claim1}. The difference is that we have to use Lemma~\ref{Lemma:basic4} instead of Lemma~\ref{Lemma:basic2}.
Furthermore, $\tilde{R}=(C^*_1,C^*_2,\ldots,C^*_{M-1},\{ k \},C^*_{M} - \{ k \})$ is path-connected only for node $k \in C^{*}_M$ such that it reaches a node in $\cup_{j=1}^{M-1} C_j^*$.  


{\bf{Proof of Uniqueness:}}

Consider two path-connected rank orderings $R=\{ C_1, C_2, \ldots, C_M \}$ and $\hat{R}=\{ \hat{C}_1, \hat{C}_2, \ldots, \hat{C}_{\hat{M}} \}$. We will prove by contradiction that $\boldsymbol{Q}$ cannot be in $D_f^c(R)$ and $D_f^c(\hat{R})$ simultaneously.

{\bf{Case I.}} There exist nodes $a$ and $b$ such that $b \prec^{R} a$ and $a \prec^{\hat{R}} b$, i.e.
\[ \left\{ \begin{array}{ll}
      b \in C^{i-1} &, \ a \in \cup_{l=i}^{M} C_{l} \\
      a \in \hat{C}^{j-1} &, \ b \in \cup_{l=j}^{\hat{M}} \hat{C}_{l}          
         \end{array} \right. . \]

If $a$ and $b$ reach a node in $C^{i-1}$ and $\hat{C}^{j-1}$ respectively, then by Lemma~\ref{Lemma:basic4}, we have
\begin{eqnarray}
\label{crit_uniquenessA1}
Q_a > Q_{C^{i-1}} \ge Q_b,
\end{eqnarray}
\begin{eqnarray}
\label{crit_uniquenessA2}
Q_b > Q_{\hat{C}^{j-1}} \ge Q_a,
\end{eqnarray}
which are contradictory. 

If $a$ does not reach a node in $C^{i-1}$, we will show that there exists a node, say $\tilde{a}$, such that it reaches a node in $C^{i-1}$ and it satisfies
\begin{eqnarray}
\label{crit_uniquenessA3}
\tilde{a} \in \cup_{l=i}^M C_{l} \ , \ \tilde{a} \in \hat{C}^{j-1} .
\end{eqnarray} 

Since $a \in \cup_{l=i}^M C_{l}$ and $a$ does not reach a node in $C^{i-1}$, all paths from node $a$ to the destination include at least a node $k \in \cup_{l=i}^M C_{l}$ such that it reaches a node in $C^{i-1}$.
On the other hand, since $a \in \hat{C}^{j-1}$ and $\hat{R}$ is path-connected, there exists a path from node $a$ to the destination which only consists of nodes from $\hat{C}^{j-1}$. Now it is clear that there exists a node, say $\tilde{a}$, on this path such that it reaches a node in $C^{i-1}$ and it satisfies (\ref{crit_uniquenessA3}).

Similarly we can show that if $b$ does not reach a node in $\hat{C}^{j-1}$, then there exists a node, say $\tilde{b}$, such that it reaches a node in $\hat{C}^{j-1}$ and it satisfies
\begin{eqnarray}
\nonumber
\tilde{b} \in C^{i-1} \ , \ \tilde{b} \in \cup_{l=j}^{\hat{M}} \hat{C}_{l} .
\end{eqnarray}    

As before we can use Lemma~\ref{Lemma:basic4} for nodes $\tilde{a}$ and $\tilde{b}$ to show that $\boldsymbol{Q}$ cannot be in $D_f^c(R)$ and $D_f^c(\hat{R})$ simultaneously.   

{\bf{Case II.}} There are no nodes $a$, $b$, such that $b \prec^{R} a$ and $a \prec^{\hat{R}} b$. 

Proof of uniqueness for this case is similar to the one provided for Lemma~\ref{existuniq}. 
We only need to show that rank orderings $R_1$, $R_2$, $\hat{R}_1$, and $\hat{R}_2$, defined identically, are path-connected.
Suppose $R_1= \{ C_1, \ldots, C_{i}, A_1, B_1, C_{i+2}, \ldots, C_M \}$ is not path-connected. Since $R= \{ C_1, \ldots, C_{i}, A_1 \cup B_1, C_{i+2}, \ldots, C_M \}$ is path-connected, the only possibility is that some node in $A_1$ has no path to the destination via nodes in $C^i$. 
On the other hand, the fact that $\hat{R}$ is path-connected implies that all nodes in $A_1$ must
have a path to the destination via nodes with rank not higher than that of $A_1$, which results in a contradiction.
In a similar way we can show that rank orderings $R_2$, $\hat{R}_1$, and $\hat{R}_2$ are path-connected.

\subsection{Proof of Lemma~\ref{contdiff}}
\label{App:contdiff}

\hspace{-0.2 in} {\bf{Lemma~\ref{contdiff}}}. {\textit{
$L^*_f(\cdot)$ is continuous and differentiable. 
}}

\begin{IEEEproof}
For all $R \in \mathcal{R}_{}$, $L_f(\cdot,R)$ is a simple quadratic function in $\boldsymbol{Q}$. Hence, to prove continuity and differentiability of $L^*_f(\cdot)$, it suffices to show that $L^*_f(\cdot)$ is continuous and differentiable 
at any $\boldsymbol{Q}$ on the hyperplane separating $D_f(R)$ and $D_f(R')$, for any adjacent rank orderings $R=( C_1, \ldots, C_i, C_{i+1}, \ldots, C_M )$ and $R'=( C_1, \ldots, C_{i-1}, C_i \cup C_{i+1},C_{i+2}, \ldots, C_M )$.

The hyperplane separating $D_f(R)$ and $D_f(R')$ is given by $\Lambda_f(\boldsymbol{Q},R,i)=\Lambda_f(\boldsymbol{Q},R',i)$. From Lemma~\ref{lesspenalty}, this hyperplane can be written as 
\begin{eqnarray}
\label{contdiff1}
f(|C^{i-1}|,|C_i|) Q_{C_i} = f(|C^{i-1}|,|C_i|+|C_{i+1}|) (Q_{C_i}+Q_{C_{i+1}}) =  f(|C^{i}|,|C_{i+1}|) Q_{C_{i+1}}.
\end{eqnarray}

On one side of this hyperplane, $L^*_f(\cdot)=L_f(\cdot,R)$, and on the other side, $L^*_f(\cdot)=L_f(\cdot,R')$.
For any $\boldsymbol{Q}$ on this hyperplane,
\begin{eqnarray}
\label{contdiff2}
\nonumber
L_f(\boldsymbol{Q},R) - L_f(\boldsymbol{Q},R') &=& f(|C^{i-1}|,|C_i|) Q^2_{C_i} + f(|C^{i}|,|C_{i+1}|) Q^2_{C_{i+1}} - f(|C^{i-1}|,|C_i|+|C_{i+1}|) (Q_{C_i}+Q_{C_{i+1}})^2 \\
\nonumber
&=& f(|C^{i-1}|,|C_i|+|C_{i+1}|) \left( (Q_{C_i}+Q_{C_{i+1}}) Q_{C_i} +  (Q_{C_i}+Q_{C_{i+1}}) Q_{C_{i+1}} - (Q_{C_i}+Q_{C_{i+1}})^2  \right) \\
&=& 0,
\end{eqnarray}
where the last equality follows from (\ref{contdiff1}).
Equation (\ref{contdiff2}) implies that $L^*_f(\cdot)$ is continuous on the hyperplane separating $D_f(R)$ and $D_f(R')$.

Similarly, to prove the differentiability of $L^*_f(\cdot)$, we have to show that $L_f(\cdot,R)$ and $L_f(\cdot,R')$ have same partial derivatives at any $\boldsymbol{Q}$ on the hyperplane separating $D_f(R)$ and $D_f(R')$.
We have,
\begin{eqnarray}
\label{contdiff8}
\frac{\partial L_f(\boldsymbol{Q},R)}{\partial Q_{k}} = 2 f(|C^{j-1}|,|C_j|) Q_{C_j} \ \ \text{for all} \ k \in C_j, \ j=1,2,\ldots,M,
\end{eqnarray}
and,
\begin{eqnarray}
\label{contdiff9}
\frac{\partial L_f(\boldsymbol{Q},R')}{\partial Q_{k}} =  \left\{ \begin{array}{ll}
      2 f(|C^{j-1}|,|C_j|) Q_{C_j} & \mbox{for all $k \in C_j, \ j \ne i,i+1$} \\
      2 f(|C^{i-1}|,|C_i|+|C_{i+1}|) (Q_{C_i}+Q_{C_{i+1}}) & \mbox{for all $k \in C_i \cup C_{i+1}$}          
         \end{array} \right. .
\end{eqnarray}

From (\ref{contdiff1}),(\ref{contdiff8}), and (\ref{contdiff9}), we have
\begin{eqnarray}
\label{contdiff10}
\nabla L_f(\boldsymbol{Q},R) = \nabla L_f(\boldsymbol{Q},R').
\end{eqnarray}

\end{IEEEproof}


\subsection{Proof of Lemmas~\ref{queuedynamic_ineq} and \ref{routemax}}
\label{App:queuedynamic_ineq}

In this appendix we prove the main steps in establishing the negative expected drift in $L_f$ under the $f$-policy. 

\hspace{-0.2 in} {\bf{Lemma~\ref{queuedynamic_ineq}}}. {\textit{
Let $R=(C_1,C_2, \ldots, C_M) \in \mathcal{R}$ and $\boldsymbol{Q}(t) \in D_f(R)$. We have
\begin{eqnarray*}
Q^2_{C_i}(t+1) - Q^2_{C_i}(t) \le \beta_f - 2 Q_{C_i}(t) (\mu^*_{C_i,out}(t) - \mu^*_{C_i,in}(t) - A_{C_i}(t)),
\end{eqnarray*}
where $\beta_f$ is a constant bounded real number.
}}

\begin{IEEEproof}

For all $C_i$, if $ Q_{C_i} \ge \frac{f(|C^{i-1}|,1)}{f(|C^{i-1}|,|C_{i}|)}$, then (\ref{Lemma:basic3}) implies that for all $k \in C_i$,
$$Q_k \ge \frac{f(|C^{i-1}|,|C_{i}|)}{f(|C^{i-1}|,1)} Q_{C_i} \ge 1.$$ 

Let $$\alpha = \max_{0 \le m < N} \max_{0 < n \le N} \frac{f(m,1)}{f(m,n)}.$$

If $ Q_{C_i}(t) \ge \alpha$, then using (\ref{QDynamic}) we obtain
\begin{eqnarray}
\label{queuedynamic01}
Q_{C_i}(t+1) \le Q_{C_i}(t) - \mu^*_{C_i,out}(t) + \mu^*_{C_i,in}(t) + A_{C_i}(t).
\end{eqnarray}
The expression above is an inequality rather than an equality because the actual number of packets routed to $C_i$ from other ranking classes may be less than $\mu^*_{C_i,in}(t)$ if there are no actual packets transmitted from the nodes in those ranking classes.

After taking the square of both sides of (\ref{queuedynamic01}) and appropriate arrangements of terms, we have
\begin{eqnarray}
\label{queuedynamic02}
\nonumber
Q^2_{C_i}(t+1) - Q^2_{C_i}(t) &\le& (\mu^*_{C_i,out}(t) - \mu^*_{C_i,in}(t) - A_{C_i}(t))^2 - 2 Q_{C_i}(t) (\mu^*_{C_i,out}(t) - \mu^*_{C_i,in}(t) - A_{C_i}(t)) \\
&\le& N^2 + N^2 (1 + A_{\max})^2 - 2 Q_{C_i}(t) (\mu^*_{C_i,out}(t) - \mu^*_{C_i,in}(t) - A_{C_i}(t)). 
\end{eqnarray}

When $ Q_{C_i}(t) < \alpha$, then again using (\ref{QDynamic}), we have
\begin{eqnarray}
\label{queuedynamic03}
Q_{C_i}(t+1) \le Q_{C_i}(t) + \mu^*_{C_i,in}(t) + A_{C_i}(t).
\end{eqnarray}

This implies that,
\begin{eqnarray}
\label{queuedynamic04}
\nonumber
Q^2_{C_i}(t+1) - Q^2_{C_i}(t) &\le& (\mu^*_{C_i,in}(t) + A_{C_i}(t))^2 + 2 Q_{C_i}(t) \mu^*_{C_i,out}(t) - 2 Q_{C_i}(t) (\mu^*_{C_i,out}(t) - \mu^*_{C_i,in}(t) - A_{C_i}(t)) \\
&\le& N^2 (1 + A_{\max})^2 + 2 \alpha_{} N - 2 Q_{C_i}(t) (\mu^*_{C_i,out}(t) - \mu^*_{C_i,in}(t) - A_{C_i}(t)). 
\end{eqnarray}

Denoting $\beta_f := N^2 + N^2 (1 + A_{\max})^2 + 2 \alpha_{} N $, (\ref{queuedynamic02}) and (\ref{queuedynamic04}) result in the assertion of the lemma.

\end{IEEEproof}




\hspace{-0.2 in} {\bf{Lemma~\ref{routemax}}}. {\textit{
Let $R=(C_1, C_2, \ldots, C_M) \in \mathcal{R}$, $\boldsymbol{Q}(t) \in D_f(R)$, and let $\{\mu^*_{ij} (t)\}_{i,j \in \Omega}$ represent routing decisions made under an $f$-policy. For any collection of routing decisions $\{\mu_{ij} (t) \}_{i,j \in \Omega}$, we have
\begin{eqnarray}
\label{App:routing01}
\sum_{i=1}^M f(|C^{i-1}|,|C_i|) Q_{C_i}(t) (\mu^*_{C_i,out}(t) - \mu^*_{C_i,in}(t)) \ge \sum_{i=1}^M f(|C^{i-1}|,|C_i|) Q_{C_i}(t) (\mu_{C_i,out}(t) - \mu_{C_i,in}(t)).
\end{eqnarray}
}}

\begin{IEEEproof}
Switching the sums in the right-hand side of (\ref{App:routing01}) and using (\ref{routingdecision}), we have 
\begin{align}
\label{routing02}
\nonumber
\lefteqn{\sum_{i=1}^M f(|C^{i-1}|,|C_i|) Q_{C_i}(t) (\mu_{C_i,out}(t) - \mu_{C_i,in}(t))} \\
\nonumber
&= \sum_{i=1}^M \sum_{k \in C_i} \sum_{j=1}^M \sum_{l \in C_j} \mu_{kl}(t) \left [ f(|C^{i-1}|,|C_{i}|) Q_{C_i}(t) - f(|C^{j-1}|,|C_{j}|) Q_{C_j}(t) \right ]\\
&\le \sum_{i=1}^M \sum_{k \in C_i} \max_{1 \le j \le M} \max_{l \in C_j} \mathbf{1}_{\{ l \in S_k(t)\}} \left [ f(|C^{i-1}|,|C_{i}|) Q_{C_i}(t) - f(|C^{j-1}|,|C_{j}|) Q_{C_j}(t) \right ].
\end{align}

Since $\boldsymbol{Q}(t) \in D_f(R)$, by Lemma~\ref{Lemma:basic1}, we have
\begin{align}
\label{routing03}
f(|C^{i-1}|,|C_{i}|) Q_{C_i}(t) \le f(|C^{i}|,|C_{i+1}|) Q_{C_{i+1}}(t) \ \ i=1,2,\ldots,M-1.
\end{align}




However, from (\ref{routing03}), the upper bound in (\ref{routing02}) is achieved under the $f$-policy, i.e.\ $\mu^*_{kl}(t)=1$ only when $l \in S_k(t)$ and $l \preceq^{R} m$ for all $m \in S_k(t)$. 

\end{IEEEproof}


\subsection{Proof of Theorem~\ref{fpolicyopt} for path-connected $f$-policy}
\label{criticalfpolicy_opt}

\hspace{-0.2 in} {\bf{Theorem~\ref{fpolicyopt}}}. {\textit{
Let $f$ be a bivariate function that satisfies conditions (C1) and (C2). Then the associated path-connected $f$-policy is throughput optimal.
}}

The proof is done in the same way as described in Subsection~\ref{Thm1Proof}. The only differences appear in the number of cones partitioning $\mathbb{R}^N_+$ and in the proof of Lemma~\ref{queuedynamic_ineq}.  

Given restrictions to path-connected rank orderings, we have the following Lyapunov function
\begin{eqnarray}
\label{critical_Lyapunov}
L^{*}_{f}(\boldsymbol{Q})=L_{f}(\boldsymbol{Q},\pi^c_f(\boldsymbol{Q}))= \sum_{R \in \mathcal{R}^{c}} L_{f}(\boldsymbol{Q},R) \mathbf{1}_{\{\boldsymbol{Q} \in D_f^c(R)\}}.
\end{eqnarray}

Proof of Lemma~\ref{queuedynamic_ineq} is modified as explained below:

If $ Q_{C_i}(t) \ge \alpha$, then for any node $k \in C_i$ that reaches a node in $C^{i-1}$, by (\ref{Lemma:basic5}), we have $Q_k \ge 1$. Now consider the nodes in $C_i$ that do not reach a node in $C^{i-1}$. Under path-connected $f$-policy, these nodes may route their packets only to nodes in $C_i$.
Hence, (\ref{queuedynamic01}) holds. The rest of the proof remains unchanged.


\subsection{Proof of Theorem~\ref{respectopt}}
\label{App:respectopt}

This short appendix establishes the following:

\hspace{-0.2 in} {\bf{Theorem~\ref{respectopt}}}. {\textit{
Suppose $\Pi_{\{ R(t) \}}$ is a priority-based routing policy that is throughput optimal. Any priority-based routing policy that respects $\Pi_{\{ R(t) \}}$ is also throughput optimal. 
}}

\begin{IEEEproof}
Suppose $\Pi'_{\{ R'(t) \}}$ is a priority-based routing policy that respects $\Pi_{\{ R(t) \}}$. 
Let $S^*_i(t) = \{k \in S_i(t): k \preceq^{R(t)} j \ \text{for all} \ j \in S_i(t) \}$ and $S^{*'}_i(t)=\{k \in S_i(t): k \preceq^{R'(t)} j \ \text{for all} \ j \in S_i(t) \}$.
Since $R'(t)$ is a refinement of $R(t)$, $S^{*'}_i(t)$ is a subset of $S^{*}_i(t)$.
By definition of the priority-based routing, $\Pi'_{\{ R'(t) \}}$ selects one of the nodes in $S^{*'}_i(t)$ as the next forwarder. 
Since $S^{*'}_i(t) \subseteq S^{*}_i(t)$, this routing decision is consistent with $\Pi_{\{ R(t) \}}$, hence, guarantees throughput optimality.
\end{IEEEproof}


\subsection{Proof of Lemma~\ref{App:Vablemma}}
\label{AppSec:Vablemma}

\begin{lemma}
\label{App:Vablemma}
For any two nodes $a$ and $b$, if $a \to b$, then 
\begin{eqnarray}
\label{App:Vab}
V_a(t) \le \frac{Q_a(t)}{p_{\min}} + V_b(t).
\end{eqnarray}
\end{lemma}

\begin{IEEEproof}
If $V_a(t) \le V_b(t)$, then (\ref{App:Vab}) follows trivially.

Now suppose $b \in U_a(t):= \{ a': a \to a' \ \text{and} \ V_{a'}(t) < V_a(t) \}$.
Without loss of generality, let $U_a(t)=\{a_1,a_2,\ldots,a_K\}$ such that $V_{a_i}(t) \le V_{a_{i+1}}(t)$ for all $i \le K$.
We can rewrite (\ref{fixpoint1}) as:
\begin{align}
\label{fixpoint2}
V_a(t) &= Q_a(t) + \sum_{i=1}^{K} V_{a_i}(t) \left(\sum_{S: i= \min \{l:\ a_l \in S \}} P(S|a)\right) + V_a(t) \left(\sum_{S: S \cap U_a(t) = \emptyset} P(S|a)\right)\\
\nonumber
&\le Q_a(t) + V_{a_1}(t) \left(\sum_{S: 1= \min \{l:\ a_l \in S \}} P(S|a)\right) + V_{a}(t) \left( \sum_{i=2}^{K} \sum_{S: i= \min \{l:\ a_l \in S \}} P(S|a) +  \sum_{S: S \cap U_a(t) = \emptyset} P(S|a)\right).
\end{align}

Rearranging terms, and putting
$$P_0=1- \left(\sum_{i=2}^{K} \sum_{S: i= \min \{l:\ a_l \in S \}} P(S|a) + \sum_{S: S \cap U_a(t) = \emptyset} P(S|a) \right) = \sum_{S: 1= \min \{l:\ a_l \in S \}} P(S|a),$$
we have,
\begin{eqnarray}
\nonumber
V_a(t) &\le& \frac{Q_a(t)}{P_0} + \frac{V_{a_1}(t)}{P_0} \left( \sum_{S: 1= \min \{l:\ a_l \in S \}} P(S|a)\right)\\
\nonumber
&=& \frac{Q_a(t)}{P_0} + V_{a_1}(t)\\
\nonumber
&\le& \frac{Q_a(t)}{p_{\min}} + V_{b}(t),
\end{eqnarray}
where the last inequality holds because $b \in U_a(t)$, and $P_0 \ge p_{\min}$.

\end{IEEEproof}

\bibliographystyle{IEEEbib}
\bibliography{fPolicy} 

\end{document}